\documentclass[11pt]{AAS}

\usepackage{bm}
\usepackage{amsmath}
\usepackage[caption=false]{subfig}
\usepackage{overcite}
\usepackage{footnpag}			      	

\usepackage{todonotes}
\usepackage{rotating}
\usepackage{enumitem}
\PaperNumber{23-232}

\title{Approximate Analytical Solutions for the Circular Restricted Three-Body Problem Including Non-Hamiltonian Solar Radiation Pressure}
\author{Hailee Hettrick\footnote{Draper Scholar, Massachusetts Institute of Technology, 70 Vassar St, Cambridge, MA 02139, (909) 838-7476, hettrick@mit.edu},
David W. Miller\footnote{Massachusetts Institute of Technology, 70 Vassar St, Cambridge, MA 02139, (617) 216-7398, millerd@mit.edu},
Begum Cannataro\footnote{The Charles Stark Draper Laboratory, 555 Technology Square, Cambridge, MA 02139, (352) 871-4053, bcannataro@draper.com}}

\begin{document}
\maketitle{}
\begin{abstract}
	The circular restricted three-body problem (CR3BP) with solar radiation pressure (SRP) has often been analyzed with assumptions made on a spacecraft's attitude, such that the problem remains Hamiltonian. These assumptions are unsatisfactorily limiting for a starshade mission since the starshade's attitude will inherently vary from the configuration that corresponds to Hamiltonian dynamics. This paper presents the derivation of the equations of motion for CR3BP with SRP that permit the application of the Lindstedt-Poincar\'e method, such that approximate solutions are produced, which may serve as invaluable trajectory design tools. Examples of periodic orbits and manifolds corresponding to three sets of attitude angles are shown and the accuracy of their seventh-order approximations is considered.
\end{abstract}

\section{Introduction}
In certain dynamical regimes of space, there exist nonlinear dynamics representing the balance of gravitational and centrifugal forces on a negligibly small body (like a spacecraft) imparted by two massive celestial bodies. The motion of these three bodies is referred to as the restricted three-body problem, and a prefix is commonly attached which describes the assumed motion of the massive bodies. Of particular interest is the circular restricted three-body problem (CR3BP). CR3BP has five equilibrium points (also known as libration or Lagrange points); these points consist of three collinear equilibrium points\cite{euler1767motu} and two off-axis equilibrium points\cite{lagrange1772essai}. These equilibrium points are typically denoted as $L_\#$ corresponding to their position.

Periodic orbits exist around these Lagrange points; however, high-precision initial conditions define these periodic orbit solutions when solved from the equations of motion\cite{farquhar1968control,howell1984three}. Therefore, a sufficiently close guess of an initial condition corresponding to a periodic orbit is required and refined by a differential corrector\cite{breakwell1979halo}. For CR3BP, this initial guess is typically solved for via Richardson's third-order approximate analytical solution of the dynamics around a specific Lagrange point\cite{richardson1980analytic,richardson1980halo}.

For trajectory design, a natural extension of understanding how to find and use periodic orbits in this regime is determining how to transfer to these orbits in a fuel-inexpensive manner. Invariant manifold theory provided the solution to the construction of transfer trajectories to periodic orbits about Lagrange points\cite{gomez1991study,gomez1991moon,gomez1993study}. Furthermore, differential correctors added to this design scheme allow for the generation of transfers from a given Earth parking orbit\cite{howell1994transfer} and from Earth-to-Halo orbit and Halo-to-Halo orbits\cite{barden1994using}. 

Given the tools described, many spacecraft missions have exploited the unique Lagrange point orbits and manifolds to achieve their scientific goals. The first spacecraft mission to utilize the periodic orbits (specifically, a Halo orbit) about a Lagrange point (Sun-Earth $L_1$ and later Sun-Earth $L_2$) was the third International Sun-Earth Explorer (ISEE-3), launched in 1978 \cite{farquhar2001flight}. The successful use of the Halo orbit in ISEE-3's mission led to the use of such periodic orbits in the vicinity of Lagrange points for several more missions: SOHO (1996) \cite{domingo1995soho}, ACE (1997) \cite{stone1998advanced}, MAP (2001) \cite{bennett2003microwave}, Genesis (2001) \cite{lo2001genesis}, Planck (2009) \cite{passvogel2010planck}, Herschel Space Observatory (2009) \cite{pilbratt2010herschel}, Chang'e 2 (2010) \cite{wu2012pre}, THEMIS-B and THEMIS-C (2011) \cite{broschart2009preliminary}, Gaia (2013) \cite{prusti2016gaia}, DSCOVR (2015) \cite{roberts2015early}, LISA Pathfinder (2015) \cite{mcnamara2008lisa}, Chang'e 4 (2018) \cite{jia2018scientific}, Chang'e 4 Relay (2018) \cite{gao2019trajectory}, and Specktr-RG (2019) \cite{kovalenko2019orbit}. Most recently, the James Webb Space Telescope (JWST) is on a Halo orbit about Sun-Earth $L_2$\cite{gardner2006james}. Orbits about Sun-Earth $L_2$ are particularly appealing for space-based telescopes due to the limited environmental disturbances (i.e. it is a highly stable thermal environment). 

During the National Academies' 2020 Decadal Survey in Astronomy and Astrophysics, two telescope missions were proposed as the successor of JWST and the Roman Space Telescope. Of interest to the research herein is the Habitable Exoplanet Observatory (HabEx) mission \cite{gaudi2020habitable}. The HabEx mission plans to use two spacecraft to characterize potentially habitable exoplanets -- one spacecraft is the space-based telescope and the second spacecraft is a large starshade. The starshade effectively acts as a solar sail due to its low mass and large area. This means that the solar radiation pressure (SRP) force exerted on the starshade is non-trivial and must be accounted for in any trajectory mission design. It is possible to treat the SRP force as a perturbation and rely on the spacecraft's control laws to account for the disturbance; however, this paper focuses on working with the perturbation rather than against it. In this manner, the work described herein can be applied to future starshade and solar sail missions.

Existing trajectory mission design tools for spacecraft in the vicinity of Lagrange points are limited in that the addition of SRP is not considered. Ultimately, these tools work with the naturally occurring dynamics of the regime and it would be beneficial to produce similar tools for the SRP problems. The existing trajectory mission design tools for CR3BP include methods to produce approximate closed-form solutions for orbits and manifolds near Lagrange points since closed-form solutions are unsolvable\cite{richardson1980halo,Jorba1998,Masdemont2005}. There exist similar approximate closed-form solutions for CR3BP with SRP; however, these solutions place assumptions on the SRP force that are unrealistic for a starshade application (e.g. assuming the surface is fixed perpendicular or parallel to the Sun-line or assuming path-independent low-thrust force)\cite{mcinnes2000strategies,Baoyin2006,Baig_2009,Lei2017}. It is important to note that for a more utilitarian approximation, no assumptions should be placed on the spacecraft attitude -- in this manner, this work considers non-Hamiltonian SRP forces.

\section{Dynamical Model}
\begin{figure}[ht!]
\centering
\includegraphics[width=0.5\textwidth]{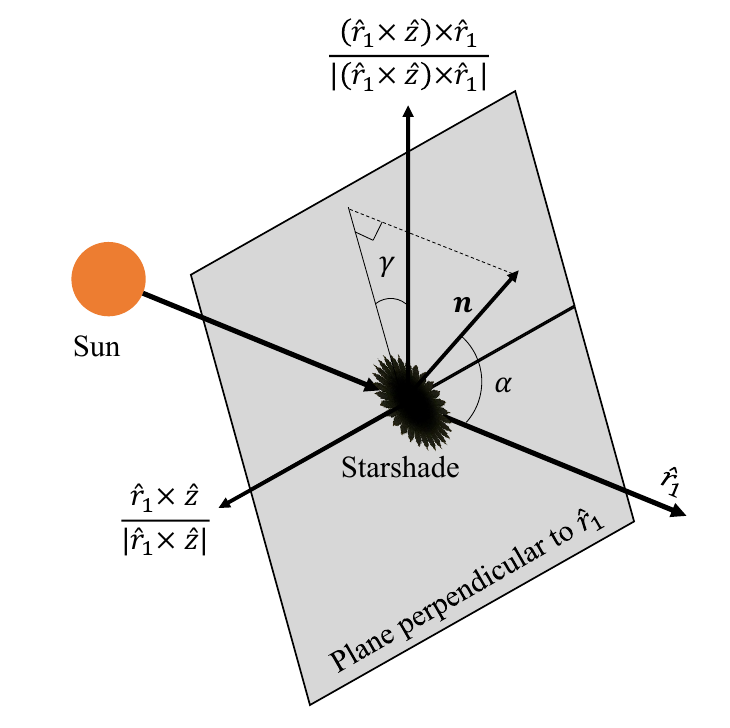}
\caption{Cone and clock angles with respect to the Sun-Starshade line}
\label{fig:sailAngles}
\end{figure}

The CR3BP with SRP equations of motion describe the motion of the third mass, which is much smaller than the two primaries, in the rotating body frame with the origin at the barycenter of the two primaries. Additionally, the third mass has a large surface that acts as a sail, experiencing a force from the solar radiation pressure. In the nominal CR3BP, the locations of the equilibrium points are a function of the mass ratio ($\mu$) of the two massive celestial bodies only. However, for CR3BP with SRP, the equilibrium points for this problem also depend on the sail's lightness number, $\beta$, and its attitude angles, $\alpha$ and $\gamma$, which are the cone and clock angles. These equilibrium points are commonly referred to as artificial equilibrium points (AEPs). The attitude angles and the starshade's normal vector $\mathbf{\hat{n}}$ are illustrated in Figure \ref{fig:sailAngles}. The domains of the attitude angles are $\alpha \in [-\pi/2,\pi/2]$ and $\gamma \in [0,\pi]$.

The CR3BP with SRP equations of motion are simply the CR3BP equations with the addition of the SRP force, denoted as $\mathbf{a}_{SRP}$:
\begin{equation}
\begin{aligned}
\ddot{x} - 2\dot{y} &= \frac{\partial \Omega}{\partial x}	 + a_{SRP,x}\\
\ddot{y} + 2\dot{x} &= \frac{\partial \Omega}{\partial y}+ a_{SRP,y}\\
\ddot{z} &= \frac{\partial \Omega}{\partial z} + a_{SRP,z}
\end{aligned}
\label{eqn:eoms}
\end{equation}
where
\begin{gather*}
\Omega = \frac{1}{2}(x^2 + y^2) + \frac{1-\mu}{r_1} + \frac{\mu}{r_2} \\
r_1 = \sqrt{(x+\mu)^2+y^2+z^2}\\
r_2 = \sqrt{x^2+y^2+z^2}	
\end{gather*}
and 
\begin{gather}
	\mathbf{a}_{SRP} = \beta \frac{1-\mu}{r_1^2} \biggr(\mathbf{\hat{r}_1} \cdot \mathbf{\hat{n}}\biggr)^2 \mathbf{\hat{n}} \\
	\mathbf{\hat{n}} = \cos (\alpha) \, \mathbf{\hat{r}_1} + \sin(\alpha) \sin (\gamma) \frac{\biggr(\mathbf{\hat{r}_1} \times \mathbf{\hat{z}}\biggr )}{|\mathbf{\hat{r}_1} \times \mathbf{\hat{z}}|} + \sin(\alpha) \cos (\gamma) \frac{\biggr( (\mathbf{\hat{r}_1} \times \mathbf{\hat{z}}) \times \mathbf{\hat{r}_1} \biggr )}{|(\mathbf{\hat{r}_1} \times \mathbf{\hat{z}}) \times \mathbf{\hat{r}_1}|}.
\end{gather}
where $\mathbf{\hat{r}_1}$ is the Sun-starshade line. 

Each element of $\mathbf{a}_{SRP}$ is
\begin{equation}
\begin{aligned}
	a_{SRP,x} &= \frac{\beta \cos^2(\alpha)(1-\mu)}{r_1^2}\biggr[\frac{(x+\mu)\cos(\alpha)}{r_1} + \frac{y \sin(\alpha)\sin(\gamma)}{r_{xy}} - \frac{z (x+\mu) \cos(\gamma)\sin(\alpha)}{r_1 r_{xy}}\biggr] \\
	a_{SRP,y} &= \frac{\beta \cos^2(\alpha)(1-\mu)}{r_1^2}\biggr [ \frac{y \cos(\alpha)}{r_1} - \frac{(x+\mu)\sin(\alpha)\sin(\gamma)}{r_{xy}} -\frac{yz\cos(\gamma)\sin(\alpha)}{r_1 r_{xy}}   \biggr ] \\
	a_{SRP,z} &= \frac{\beta \cos^2(\alpha)(1-\mu)}{r_1^3} \biggr [ z \cos(\alpha) + \frac{((x+\mu)^2+y^2) \cos(\gamma)\sin(\alpha)}{r_{xy}} \biggr ]
\end{aligned}
\label{eqn:asrp}
\end{equation}
where $r_{xy} = \sqrt{(x+\mu)^2+y^2}$. The right-hand side of the equations of motion can be denoted as $\mathbf{F}=\nabla\Omega + \mathbf{a}_{SRP}$. In CR3BP, $\Omega$ is typically referred to as the effective potential. Therefore, for CR3BP with non-Hamiltonian SRP, the relevant effective potential can be denoted as $\Omega^* = \Omega + \mathbf{a}_{SRP} \cdot [x, \, y, \, z]^T$.

The derivation of the CR3BP equations can be found in several texts, but Koon et al.\cite{koon2000dynamical} is an excellent source.  Multiple texts also derive the equations of motion (see Eq.~\eqref{eqn:eoms}) for the SRP case; however, C. McInnes\cite{mcinnes1994solar} and A.I.S. McInnes\cite{mcinnes2000strategies} were mostly used for this work.
\subsection{Artificial Equilibrium Points}

Due to the addition of the SRP force, new equilibrium solutions now exist that differ from the Lagrange points of the circular restricted three-body problem. These artificial equilibrium points vary with the values of $\mu$, $\beta$, $\alpha$, and $\gamma$. Parameterization of the equilibrium solutions into surfaces is illustrated by A.I.S. McInnes\cite{mcinnes2000strategies}. Since the starshade has a known $\beta$ value, evaluating how the AEPs vary with $\alpha$ and $\gamma$ only is of consequence for this work. 
\begin{figure}[!htp]
\centering
{%
  \includegraphics[width=45mm]{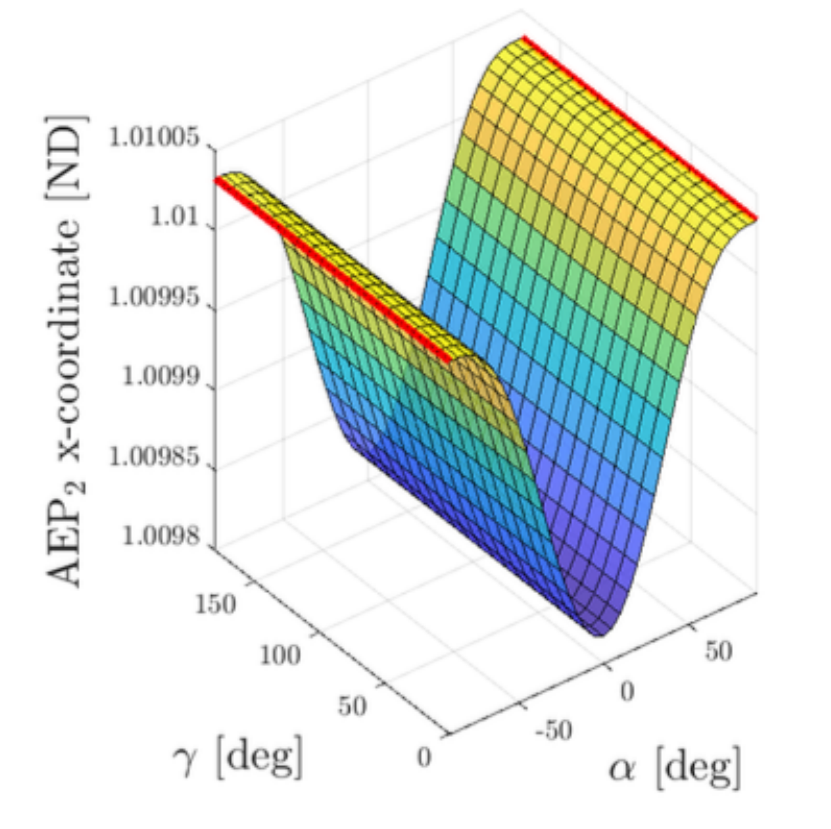}%
}\quad
{%
  \includegraphics[width=45mm]{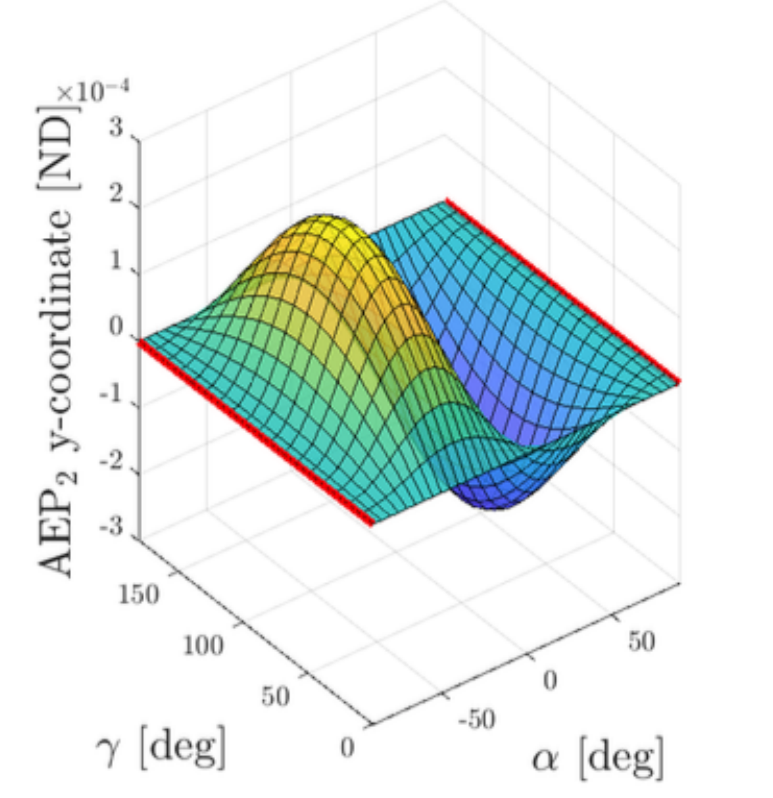}
}\quad
{%
  \includegraphics[width=45mm]{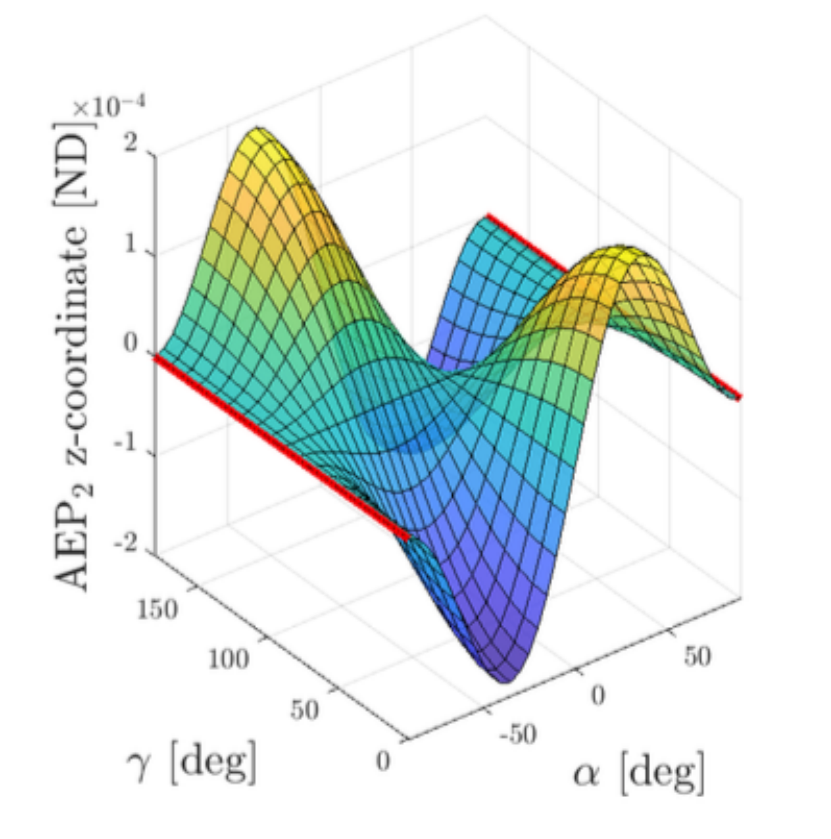}
  }
\caption{Position of $H_2$ as it varies with respect to $\alpha$ and $\gamma$ for a constant $\beta$}
\label{fig:l2Location}
\end{figure}
Equilibrium solutions correspond to positions with zero velocity and acceleration in the rotating reference frame. Therefore, equilibria are computed from the relationship
\begin{equation}
- \mathbf{\nabla} \Omega = \beta \frac{1-\mu}{r_1^2} \biggr(\hat{r}_1 \cdot \mathbf{\hat{n}}\biggr)^2 \mathbf{\hat{n}}.
\end{equation}
Newton's method may be used to determine the artificial equilibrium points. An initial guess must be provided when using Newton's method, so the corresponding Lagrange point position from CR3BP is used. An example of the possible positions for $H_2$ (analogous to $L_2$ in CR3BP) in the Sun-Earth system as a function of varying attitude angles for a constant sail lightness number $\beta$ is shown in Figure \ref{fig:l2Location}.

\section{Finding Equivalent Formulation of the Equations of Motion}
\begin{figure}[h]
\includegraphics[width=0.7\textwidth]{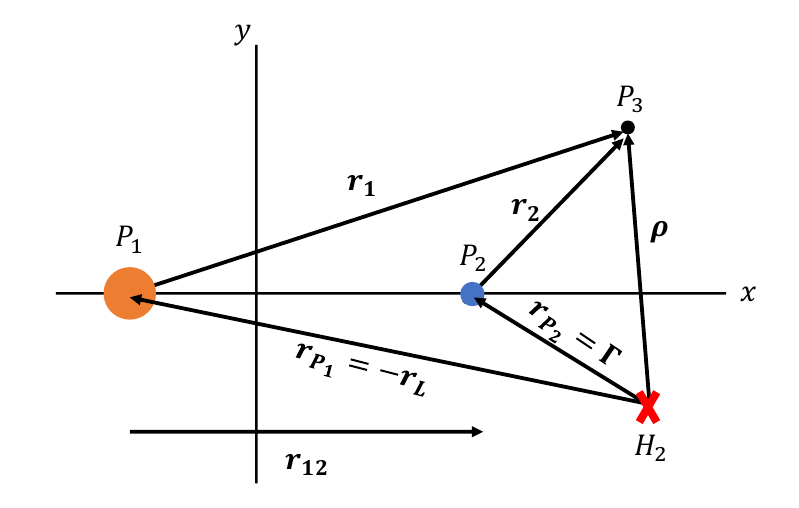}
\centering
\caption{Geometry of system with respect to $H_2$}
\label{fig:alphaFree}
\end{figure}

To find approximate closed-form solutions for CR3BP with SRP, a different but equivalent formulation of the equations of motion is sought. First, the dynamics of CR3BP with SRP is translated from being about the barycenter of the two massive celestial bodies to about the AEP of interest. The position vector $\boldsymbol{\rho} = [x, \, y, \, z]^T$ describes the position of the third body with respect to the AEP. In other words, the system is translated about the AEP located at $[H_x, H_y, H_z]^T$ with respect to the barycenter. The geometry of the problem is illustrated in Figure \ref{fig:alphaFree}.
\subsection{Lagrangian}
Next, it is necessary to find the Lagrangian of the system as the difference between the kinetic and potential energies of the system, yielding
\begin{gather}
\begin{aligned}
	\mathcal{L} &= \frac{1}{2}(\boldsymbol{\dot{\rho}} \cdot \boldsymbol{\dot{\rho}})
	+ (1-\mu)\biggr(\frac{1}{|\boldsymbol{\rho}-\mathbf{r_{P_1}}|} -\frac{\mathbf{r_{P_1}} \cdot \boldsymbol{\rho}}{r_{P_1}^3}  \biggr) +
	\mu \biggr[ \frac{1}{|\mathbf{r_{P_2}} - \boldsymbol{\rho}|} -\frac{\mathbf{r_{P_2}} \cdot \boldsymbol{\rho}}{r_{P_2}^3}\biggr] \\
	&- (1-\mu) \beta \cos^2(\alpha) \biggr(\frac{\mathbf{\hat{n}(r_L) }\cdot \boldsymbol{\rho}}{r_{P_1}^2}    \biggr)
\end{aligned}
\label{eqn:lp1}
\end{gather}
where $\mathbf{r_{P_1}} = [-H_x-\mu, \, -H_y, \, -H_z]^T$,  $\mathbf{r_{P_2}} = \mathbf{\Gamma} = [1-H_x-\mu, \, -H_y, \, -H_z]^T$, and $\mathbf{\hat{n}(r_L)}$ is the evaluation of the normal vector at $\mathbf{r_L} = [H_x + \mu, \, H_y, \, H_z]^T$. Refer to Figure \ref{fig:alphaFree} for an illustration of how these vectors relate.

Richardson provides some details on the Lagrangian formulation for CR3BP\cite{Richardson1980}; however, the Lagrangian here differs since the equilibrium points for CR3BP with SRP are no longer (generally) collinear with the $x$-axis. The result of the difference is the final term in the Lagrangian in Eq.~\eqref{eqn:lp1}. Note that when the starshade is edge-on to the Sun-line, the Lagrangian reduces to that which is characteristic of CR3BP ($\alpha=\pm 90^\circ$ and $\gamma = 0^\circ$ when edge-on and the SRP force is nulled). Furthermore, note that when the starshade is perpendicular to the Sun-line ($\alpha = 0^\circ$), the Lagrangian found by A.I.S. McInnes\cite{mcinnes2000strategies} is recreated.

\subsection{Legendre Polynomials}
Prior to applying Lagrange's equation to the Lagrangian in Eq.~\eqref{eqn:lp1}, a third-body perturbation simplification\cite{brouwer2013methods} is applied via the Legendre polynomial identity
\begin{equation}
\frac{1}{\sqrt{(x-A)^2 + (y-B)^2+ (z-C)^2}} = \frac{1}{D} \sum_{n=0}^{\infty} \biggr( \frac{\rho}{D} \biggr)^n P_n \biggr( \frac{Ax+By+Cz}{D\rho}\biggr)
\label{eqn:lpIdentity}
\end{equation}
where $D^2 = A^2+B^2+C^2$ and $\rho^2 = x^2 + y^2 + z^2$. For this identity to be used, the coordinates must be scaled to satisfy the numerical properties of Legendre polynomials (i.e. $\frac{\rho}{D} < 1$ and $\biggr(\frac{Ax+By+Cz}{D\rho}\biggr) < 1$). Treating $x,y,z$ as having already been scaled by $|\mathbf{\Gamma}|$, the other lengths scaled are the elements of the position vectors of $P_1$ and $P_2$ with respect to the AEP.
\begin{equation}
	A_1 = -\frac{H_x + \mu}{|\mathbf{\Gamma}|}, \quad A_2= \frac{1-H_x - \mu}{|\mathbf{\Gamma}|} \quad B = -\frac{H_y}{|\mathbf{\Gamma}|}, \quad C = - \frac{H_z}{|\mathbf{\Gamma}|}.
	\label{eqn:lengthsScaled}
\end{equation}
Applying the Legendre polynomial identity (Eq.~\eqref{eqn:lpIdentity}) to the two gravitational potential terms in Eq.~\eqref{eqn:lp1} yields
\begin{equation}
	\frac{1}{|\boldsymbol{\rho}-\mathbf{r_{P_1}}|} -\frac{\mathbf{r_{P_1}} \cdot \boldsymbol{\rho}}{r_{P_1}^3} = \frac{1}{D_1} + \sum_{n=2}^\infty \biggr( \frac{\rho^n}{D_1^{n+1}}\biggr) P_n\biggr(\frac{A_1x+By+Cz}{D_1\rho}\biggr)
\end{equation}
\begin{equation}
	\frac{1}{|\boldsymbol{\rho}-\mathbf{r_{P_2}}|} -\frac{\mathbf{r_{P_2}} \cdot \boldsymbol{\rho}}{r_{P_2}^3} = \frac{1}{D_2} + \sum_{n=2}^\infty \biggr( \frac{\rho^n}{D_2^{n+1}}\biggr) P_n\biggr(\frac{A_2x+By+Cz}{D_2\rho}\biggr).
\end{equation}
With these expressions, a new and equivalent formulation of the Lagrangian is 
\begin{equation}
\begin{aligned}
	\mathcal{L} &= \frac{1}{2}(\boldsymbol{\dot{\rho}} \cdot \boldsymbol{\dot{\rho}})
	+ (1-\mu)\biggr(\frac{1}{D_1} + \sum_{n=2}^\infty \biggr( \frac{\rho^n}{D_1^{n+1}}\biggr) P_n\biggr(\frac{A_1x+By+Cz}{D_1\rho}\biggr)
\biggr) + \\
	&\mu \biggr[ \frac{1}{D_2} + \sum_{n=2}^\infty \biggr( \frac{\rho^n}{D_2^{n+1}}\biggr) P_n\biggr(\frac{A_2x+By+Cz}{D_2\rho}\biggr)\biggr]
	- (1-\mu) \beta \cos^2(\alpha) \biggr(\frac{\mathbf{\hat{n}(r_L) }\cdot \boldsymbol{\rho}}{r_{P_1}^2}    \biggr).
\end{aligned}
\label{eqn:lag2}
\end{equation}

\subsection{Applying Lagrange's Equations}
Now, the Lagrangian in Eq.~\eqref{eqn:lag2} may be used to find the equations of motion. Note that the non-Hamiltonian SRP force means the problem is non-conservative. Thus, when applying Lagrange's equations, generalized forces are used to account for the contributions from SRP. Recall that Lagrange's equations with generalized forces are
\begin{equation}
	\frac{d}{dt}\frac{\partial \mathcal{L}}{\partial \dot{q}_j} - \frac{\partial \mathcal{L}}{\partial q_j} = Q_j, \quad j=x,y,z
\end{equation}
where $q_j$ represents the variables of the problem -- in this problem, we apply the Lagrange equations to $x$, $y$, and $z$. Finding the generalized forces is straightforward for this problem via the following equation
\begin{equation*}
Q_j = \sum \mathbf{a}_{SRP} \cdot \biggr(\frac{\partial \boldsymbol{\rho}}{\partial q_j}\biggr)	, \quad j=x,y,z.
\end{equation*}
Recall that $\mathbf{a}_{SRP} = \beta \frac{1-\mu}{r_1^2}\cos^2(\alpha) \mathbf{\hat{n}}$. For each $q_j$ variable, the following is calculated
\begin{gather}
Q_x = \sum \mathbf{a}_{SRP} \cdot \biggr(\frac{\partial \boldsymbol{\rho}}{\partial x}\biggr) = a_{SRP,x} \\
Q_y = \sum \mathbf{a}_{SRP} \cdot \biggr(\frac{\partial \boldsymbol{\rho}}{\partial y}\biggr) = a_{SRP,y} \\
Q_z = \sum \mathbf{a}_{SRP} \cdot \biggr(\frac{\partial \boldsymbol{\rho}}{\partial z}\biggr) = a_{SRP,z}
\end{gather}
where the subscripts $x$, $y$, and $z$ correspond to the first, second, and third elements of $\mathbf{a}_{SRP}$. Since $\boldsymbol{\rho}$ is in a rotating frame, differentiation occurs via the transport theorem.
\begin{gather*}
	\frac{d}{dt}\boldsymbol{\rho} = \biggr[ \biggr(\frac{d}{dt}\biggr)_r + \boldsymbol{\omega} \times \biggr ] \boldsymbol{\rho}
\end{gather*}
In this problem, the angular velocity is nondimensionalized to unity: $\boldsymbol{\omega} = [0,\, 0,\, 1]^T$. Therefore, the following is the expression for the first and second-time derivatives of $\boldsymbol{\rho}$.
\begin{gather*}
	\frac{d}{dt}\boldsymbol{\rho}	=\left [
	\begin{array}{ccc}
		\dot{x}-y \\
		\dot{y}+x \\
		\dot{z}
	\end{array}\right ] \\ 
	\frac{d^2}{dt^2}\boldsymbol{\rho} =\left [
	\begin{array}{ccc}
		\ddot{x}-2\dot{y} -x\\
		\ddot{y}+2\dot{x}-y \\
		\ddot{z}
	\end{array}\right ]
\end{gather*}
The partial derivative of the final term $(1-\mu) \beta \cos^2(\alpha) \biggr(\frac{\mathbf{\hat{n}(r_L) }\cdot \boldsymbol{\rho}}{r_{P_1}^2}    \biggr)$ with respect to each coordinate $q_j$ is
\begin{equation}
\begin{aligned}
	\frac{\partial G}{\partial x} &=
	 \biggr[\frac{\beta (-1 + \mu) \cos(\alpha)^2}{r_L^2} \biggr(
	 \frac{(H_x + \mu) \cos(\alpha)}{r_L} - \frac{H_z (H_x + \mu) \frac{\sqrt{H_y^2 + (H_x + \mu)^2}}{r_L} \cos(\gamma) \sin(\alpha)}{H_y^2 + (H_x + \mu)^2} + \\
	 &\frac{H_y \sin(\alpha) \sin(\gamma)}{\sqrt{H_y^2 + (H_x + \mu)^2}}
	 \biggr) 
	 \biggr] \\
      \frac{\partial G}{\partial y} &=
       \biggr[\frac{\beta (-1 + \mu) \cos(\alpha)^2}{ r_L^2} \biggr( \frac{H_y \cos(\alpha)}{r_L} - 
       \frac{H_y H_z \frac{\sqrt{H_y^2 + (H_x + \mu)^2}}{r_L}\cos(\gamma) \sin(\alpha)}{H_y^2 + (H_x + \mu)^2} - \\
       &\frac{(H_x + \mu) \sin(\alpha) \sin(\gamma)}{\sqrt{H_y^2 + (H_x + \mu)^2}} \biggr) \biggr] \\
      \frac{\partial G}{\partial z} &=
      \biggr[ \frac{\beta (-1 + \mu) \cos(\alpha)^2}{r_L^3}\biggr(H_z \cos(\alpha) + \sqrt{H_y^2 + (H_x + \mu)^2} \cos(\gamma) \sin(\alpha) \biggr) \biggr]
\end{aligned}
\end{equation}
where
\begin{equation*}
	r_L = \sqrt{(H_x+\mu)^2 + H_y^2 + H_z^2}.
\end{equation*}
\subsection{Recurrence Relations}
While it is possible to find the partial derivatives of the Legendre polynomial power series directly, it better suits the purpose (and computational needs) of finding approximate analytical solutions to define a recurrence relation for the power series and find its partial derivative. The Legendre polynomial recurrence property
\begin{equation}
	(n+1) P_{n+1}(S) = (2n+1) S P_n(S) - n P_{n-1}(S)
\end{equation}
is applied to $T_n = (\frac{\rho}{D})^n P_n\left( \frac{Ax+By+Cz}{D\rho} \right)$. This yields the following recurrence relation for $T_n$.
\begin{equation}
	T_n =  \frac{(2n-1)}{n} \left(\frac{Ax+By+Cz}{D^2} \right) T_{n-1} - \frac{(n-1)}{n} \left(\frac{\rho}{D}\right)^2 T_{n-2}
\end{equation}
This relation can be initialized with $T_0 = 1$ and $T_1 = \frac{Ax+By+Cz}{D^2}$. The recurrence relation for its partial derivative with respect to $q_j$ ($R_n^{q_j} = \frac{\partial T_n}{\partial q_j}$) is
\begin{equation}
	\begin{aligned}
	R_n^{q_j} &=  \frac{(2n-1)}{n}\biggr\{ \biggr(\frac{Q_{q_j}}{D^2}  \biggr) T_{n-1}\biggr(\frac{Ax+By+Cz}{D\rho}\biggr) +   \biggr(\frac{Ax+By+Cz}{D\rho}\biggr)R_{n-1}^{q_j}\biggr\} - \\
	&\frac{(n-1)}{n}\biggr\{\frac{2q_j}{D^2}  T_{n-2}\biggr(\frac{Ax+By+Cz}{D\rho}\biggr) + \biggr(\frac{\rho}{D}\biggr)^2 R_{n-2}^{q_j} \biggr\}
\end{aligned}
\end{equation}
where $Q_{q_j}$ is $A$ when $q_j=x$, $B$ when $q_j=y$, and $C$ when $q_j = z$. This relation is initialized with $R_1^{q_j} = \frac{Q_{q_j}}{D^2} $ and $R_2^{q_j} = \frac{3Q_{q_j}}{D^4}(Ax+By+Cz) -\frac{q_j}{D^2}$.

Therefore, the following equivalency has been determined for the Legendre polynomial power series for the gravitational potential terms, which offers computational cost savings.
\begin{equation}
\begin{aligned}
	\frac{\partial}{\partial q_j} \biggr( (1-\mu) \sum_{n=2}^\infty \frac{\rho^n}{D_1^{n+1}} P_n \biggr( \frac{A_1x+By+Cz}{D_1 \rho} \biggr) + \mu \sum_{n=2}^\infty \frac{\rho^n}{D_2^{n+1}} P_n \biggr( \frac{A_2x+By+Cz}{D_2 \rho} \biggr) \biggr) \\ 
	= \sum_{n=2}^\infty \biggr[ \frac{1-\mu}{D_1} R_n^{q_j,1} + \frac{\mu}{D_2} R_n^{q_j,2} \biggr] 
	\end{aligned}
\end{equation}

\subsection{SRP Contributions}
The Legendre polynomial identity in Eq.~\eqref{eqn:lpIdentity} is also applied to the SRP terms: $a_{SRP,x}$, $a_{SRP,y}$, and $a_{SRP,z}$. The translation and scaling previously mentioned are first applied to the terms shown in Eq.~\eqref{eqn:asrp}, making use of the constants defined in Eq.~\eqref{eqn:lengthsScaled}. Denoting $f = \beta \cos^2 \alpha (1-\mu)$, the SRP terms are now expressed as
\begin{equation}
\begin{aligned}
	a_{SRP,x} &= \frac{f}{r_1^2}\biggr( \frac{\cos \alpha (x -A_1)}{r_1 }+ \frac{(y-B) \sin \alpha \sin\gamma}{r_{xy}} - \frac{(z-C) \cos \gamma \sin\alpha (x -A_1)}{r_1 r_{xy}} \biggr) \\
	a_{SRP,y} &= \frac{f}{r_1^2}\biggr( \frac{ (y-B) \cos\alpha}{r_1} - \frac{(x -A_1)\sin\alpha \sin\gamma }{r_{xy}} - \frac{(y-B) (z-C) \cos \gamma \sin\alpha}{r_1 r_{xy} }\biggr)\\
a_{SRP,z}&= \frac{f}{r_1^3}\biggr( (z-C) \cos\alpha + \frac{((x -A_1)^2 + (y-B)^2)\cos \gamma \sin\alpha }{r_{xy} }\biggr).
\end{aligned}
\end{equation}
Each instance of $\frac{1}{r_1}$ and $\frac{1}{r_{xy}}$ is equivalent to the following series.
\begin{gather*}
\begin{aligned}
	\frac{1}{r_1} = \frac{1}{\sqrt{ (x-A_1)^2 + (y-B)^2 + (z-C)^2}} &= \frac{1}{D} \sum_{n=0}^\infty \biggr( \frac{\rho}{D}\biggr)^n P_n\biggr(\frac{A_1x+By+Cz}{D\rho}\biggr) \\
	\frac{1}{r_1}&=\sum_{n=0}^\infty \biggr( \frac{\rho^n}{D^{n+1}}\biggr) P_n\biggr(\frac{A_1x+By+Cz}{D\rho}\biggr) 
	\end{aligned} \\
\begin{aligned}
	\frac{1}{r_{xy}} = \frac{1}{\sqrt{ (x-A_1)^2 + (y-B)^2}} &= \frac{1}{D_{xy}} \sum_{n=0}^\infty \biggr( \frac{\rho_{xy}}{D_{xy}}\biggr)^n P_n\biggr(\frac{A_1x+By}{D_{xy}\rho_{xy}}\biggr) \\
	\frac{1}{r_{xy}} &= \sum_{n=0}^\infty \biggr( \frac{\rho_{xy}^n}{D_{xy}^{n+1}}\biggr) P_n\biggr(\frac{A_1x+By}{D_{xy}\rho_{xy}}\biggr) 
\end{aligned}
\end{gather*}
where $\rho_{xy}^2 = x^2 + y^2$ and $D_{xy}^2 = A_1^2+B^2$. Utilizing the recurrence relation $T_n$ found for the gravitational potential terms, it is clear that $\frac{1}{r_1} = \sum_{n=0}^\infty T_n$, and a similar recurrence relation can be found $\frac{1}{r_{xy}} = \sum_{n=0}^\infty T_{xy,n}$. Thus, the SRP terms are  equivalently expressed as
\begin{equation}
\begin{aligned}
	a_{SRP,x} &= f\biggr[ (x-A_1) \cos(\alpha) \left ( \sum_{n=0}^\infty T_n^3 \right) + (y-B)\sin(\alpha)\sin(\gamma) \left( \sum_{n=0}^\infty T_n^2 \right) \left( \sum_{n=0}^\infty T_{xy,n} \right) - \\
	&(z-C)(x-A_1) \cos(\gamma)\sin(\alpha) \left ( \sum_{n=0}^\infty T_n^3 \right)\left( \sum_{n=0}^\infty T_{xy,n} \right)\biggr] \\
	a_{SRP,y} &= f\biggr[ (y-B) \cos(\alpha) \left ( \sum_{n=0}^\infty T_n^3 \right) - (x-A_1)\sin(\alpha)\sin(\gamma) \left( \sum_{n=0}^\infty T_n^2 \right) \left( \sum_{n=0}^\infty T_{xy,n} \right) - \\
	&(z-C)(y-B) \cos(\gamma)\sin(\alpha) \left ( \sum_{n=0}^\infty T_n^3 \right)\left( \sum_{n=0}^\infty T_{xy,n} \right)\biggr] \\
	a_{SRP,z} &= f\biggr[ (z-C) \cos(\alpha) \left ( \sum_{n=0}^\infty T_n^3 \right) + ((x -A_1)^2 + (y-B)^2) \cos(\gamma)\sin(\alpha) \left ( \sum_{n=0}^\infty T_n^3 \right)\left( \sum_{n=0}^\infty T_{xy,n} \right)\biggr]
\end{aligned}
\end{equation}

\subsection{Equations of Motion for Lindstedt-Poincar\'e Method}
Prior to constructing the equations of motion to be used by the Lindstedt-Poincar\'e method, some terms are extracted from the power series to better facilitate the iterative nature of the method. The second partial derivatives of the effective potential $\Omega^*$, denoted by $\Omega^*_{ij}$, are equivalent to the scaled sum of the first-order terms from the gravitational potential power series and from the SRP power series. For example,
\begin{equation}
	\Omega_{xx}^* = \frac{1}{|\mathbf{\Gamma}|^3}  \frac{\partial}{\partial x} \biggr[ \left(a_{SRP,x}+ \frac{1-\mu}{D_1}R_2^{x,1}+ \frac{\mu}{D_2} R_2^{x,2} \right)\biggr|_{ord=1} \biggr]
\end{equation}
where the first-order terms of the $R_n$ recurrence relation are explicitly equivalent to the values evaluated at $n=2$.

Moreover, when the SRP power series are evaluated at $n=0$, there exist constant values
\begin{equation}
\begin{aligned}
	a_{SRP,x} \biggr|_{ord=0} &= \frac{f}{D_1^2}\biggr( -\frac{A_1 \cos \alpha}{D_1 } - \frac{B \sin \alpha \sin\gamma}{D_{xy}} - \frac{A_1 C \cos \gamma \sin\alpha}{D_1 D_{xy}} \biggr) \\
	a_{SRP,y}\biggr|_{ord=0} &= \frac{f}{D_1^2}\biggr( -\frac{ B \cos\alpha}{D_1} + \frac{A_1\sin\alpha \sin\gamma }{D_{xy}} - \frac{B C \cos \gamma \sin\alpha}{D_1 D_{xy} }\biggr)\\
a_{SRP,z}\biggr|_{ord=0}&= \frac{f}{D_1^3}\biggr( -C \cos\alpha + \frac{(A_1^2 + B^2)\cos \gamma \sin\alpha }{D_{xy} }\biggr)
\end{aligned}
\end{equation}
which are combined with the constants found from applying Lagrange's equations (the $\frac{\partial G}{\partial q_j}$ terms). Ultimately, the consequence of the AEP \textit{not} being guaranteed to be collinear results in a particular solution element in the linearized equations of motion (shown in the next section), which does not arise in the CR3BP case.

Now, the equations of motion for the Lindstedt-Poincar\'e method are
\begin{equation}
\begin{aligned}
	\ddot{x}-2\dot{y}  -\Omega_{xx}^* x - \Omega_{xy}^* y - \Omega_{xz}^* z &= \frac{1}{|\mathbf{\Gamma}|^3} \biggr( \frac{\partial G}{\partial x} + a_{SRP,x}\biggr|_{ord=0} + \\
	 &a_{SRP,x}\biggr|_{ord\geq2}  +  \sum_{n\geq3} \biggr[ \frac{1-\mu}{D_1} R_n^{x,1} + \frac{\mu}{D_2} R_n^{x,2} \biggr] \biggr) \\
	\ddot{y}+2\dot{x} -\Omega_{yx}^* x - \Omega_{yy}^* y - \Omega_{yz}^* z &= \frac{1}{|\mathbf{\Gamma}|^3} \biggr( \frac{\partial G}{\partial y} +  a_{SRP,y}\biggr|_{ord=0}+ \\
	 &a_{SRP,y}\biggr|_{ord\geq2}  +  \sum_{n\geq3} \biggr[ \frac{1-\mu}{D_1} R_n^{y,1} + \frac{\mu}{D_2} R_n^{y,2} \biggr] \biggr) \\
	\ddot{z}-\Omega_{zx}^* x - \Omega_{zy}^* y - \Omega_{zz}^* z &= \frac{1}{|\mathbf{\Gamma}|^3} \biggr(  \frac{\partial G}{\partial z} + a_{SRP,z}\biggr|_{ord=0}  + \\
	&a_{SRP,z}\biggr|_{ord\geq2}  + \sum_{n\geq3} \biggr[ \frac{1-\mu}{D_1} R_n^{z,1} + \frac{\mu}{D_2} R_n^{z,2} \biggr] \biggr).
\end{aligned}
\label{eqn:theBigEom}	
\end{equation}

\subsection{Linearization}
The linearized equations of motion are found by ignoring the higher-order terms in Eq.~\eqref{eqn:theBigEom}, such that
\begin{equation}
	\begin{aligned}
		&\ddot{x} -2\dot{y} -\Omega_{xx}^* x - \Omega_{xy}^* y - \Omega_{xz}^* z = \frac{1}{|\mathbf{\Gamma}|^3} \biggr( \frac{\partial G}{\partial x} + a_{SRP,x}\biggr|_{ord=0} \biggr) \\
		&\ddot{y} +2\dot{x} -\Omega_{yx}^* x - \Omega_{yy}^* y - \Omega_{yz}^* z = \frac{1}{|\mathbf{\Gamma}|^3} \biggr( \frac{\partial G}{\partial y} + a_{SRP,y}\biggr|_{ord=0} \biggr) \\
		&\ddot{z} -\Omega_{zx}^* x - \Omega_{zy}^* y - \Omega_{zz}^* z = \frac{1}{|\mathbf{\Gamma}|^3} \biggr( \frac{\partial G}{\partial z} + a_{SRP,z}\biggr|_{ord=0} \biggr)
	\end{aligned}
\end{equation}
 
The characteristic equation of the linearized equations of motion is of the form $P(\lambda) = \lambda^6  + a_1 \lambda^4 + a_2 \lambda^3 + a_3 \lambda^2 + a_4 \lambda + a_5 = 0$, which yields six eigenvalues of the following characteristics:
\begin{enumerate}[noitemsep,topsep=0pt]
	\item Two real eigenvalues of the same magnitude and opposite sign: $(\lambda_r, \, -\lambda_r)$
	\item Two pairs of complex conjugate pairs of eigenvalues: $(\omega_r + i \omega_0, \omega_r - i \omega_0)$ and $(v_r+ i v_0, v_r - i v_0)$. These pairs become purely imaginary when $\alpha = 0$, $\alpha = \pm\pi/2$, $\gamma = 0$, or $\gamma = \pi$.
\end{enumerate}
The eigenvectors associated with the complex pairs are denoted by $\mathbf{u}_1 \pm i \mathbf{w}_1$ and $\mathbf{u}_2 \pm i \mathbf{w}_2$, whereas the eigenvectors associated with the real eigenvalues are denoted by $\mathbf{v}_1$ and $\mathbf{v}_2$. Therefore, the general solution of the linearized equations of motion around the AEP is of the form:
\begin{equation}
\begin{aligned}
	x(t,\alpha_1,\alpha_2,\alpha_3,\alpha_4) &= \alpha_1 e^{\lambda_r t} + \alpha_2 e^{-\lambda_r t} + \alpha_3 e^{\omega_r t}\biggr(\cos (\omega_0 t + \phi_1) + k_{14} \sin(\omega_0 t + \phi_1) \biggr) + \\ 
	&\alpha_4 e^{v_r t} \biggr( k_5 \cos(v_0 t + \phi_2) + k_6 \sin(v_0 t + \phi_2) \biggr) + k_{15} \\
	y(t,\alpha_1,\alpha_2,\alpha_3,\alpha_4) &= k_1 \alpha_1 e^{\lambda_r t} + k_2 \alpha_2 e^{-\lambda_r t} +  \alpha_3  e^{\omega_r t}\biggr( k_7 \cos (\omega_0 t + \phi_1) + k_8 \sin(\omega_0 t + \phi_1) \biggr) + \\
	&\alpha_4 e^{v_r t} \biggr( k_9 \cos(v_0 t + \phi_2) + k_{10} \sin(v_0 t + \phi_2) \biggr) + k_{16}\\
	z(t,\alpha_1,\alpha_2,\alpha_3,\alpha_4) &= k_3 \alpha_1 e^{\lambda_r t} + k_4 \alpha_2 e^{-\lambda_r t} +  \alpha_3  e^{\omega_r t}\biggr( k_{11} \cos (\omega_0 t + \phi_1) + k_{12} \sin(\omega_0 t + \phi_1) \biggr) + \\
	& \alpha_4 e^{v_r t} \biggr( \cos(v_0 t + \phi_2) + k_{13} \sin(v_0 t + \phi_2) \biggr) + k_{17}
\end{aligned}	
\label{eqn:lin}
\end{equation}
where the amplitudes $\alpha_1$, $\alpha_2$, $\alpha_3$, and $\alpha_4$ are user-defined and the $k_\#$ coefficients are functions of $\Omega_{ij}^*$, $\lambda_0$, $\nu_0$, and $\omega_0$ and the eigenvectors.

Note that there is a reduction that has been made to these linearized equations of motion to guarantee it is solvable for all non-Hamiltonian SRP forces. This reduction is the assumption that the real part of the eigenvalues corresponding to $\omega_0$ and $\nu_0$ ($\omega_r$ and $\nu_r$) is exactly zero. The justification for this assumption is an observation from solving for the eigenvalues of many attitude angle and sail lightness number combinations and observing the real part is always close to zero\footnotemark\footnotetext{A numerical study was executed which solved for the frequencies for a mesh grid of the ranges of $\alpha$ and $\gamma$. The absolute average of all calculated $\omega_r$ values was $1 \times 10^{-5}$; the absolute average of all calculate $\nu_r$ values was $1 \times 10^{-7}$.}. Given that these real values would be the argument in an exponential in the solution (namely, $e^{\omega_r t}$ and $e^{\nu_r t}$), the fact the values are always near zero indicates these terms are equivalent to a value of one. In this manner, this reduction serves as an additional approximation.

\section{Methodology for Finding Approximate Solutions}
Given the CR3BP with SRP equations of motion and the generalized solution of the linearized equations of motion, perturbation methods can now be used to determine approximate analytical closed-form solutions to CR3BP with a non-Hamiltonian SRP force. Similar to Richardson\cite{richardson1980analytic}, Lei\cite{Lei2017}, and Masdemont \cite{Masdemont2005}, the Lindstedt-Poincar\'e method is used to find higher-order solutions.

The approximate analytical solutions for periodic orbits and their manifolds can be found in the following form.
\begin{equation}
\begin{aligned}
	x(t,\alpha_1,\alpha_2,\alpha_3,\alpha_4) = \sum e^{(i-j)\theta_3} \biggr(x_{ijkm}^{pq} \cos(p\theta_1 + q\theta_2) + \bar{x}_{ijkm}^{pq} \sin(p\theta_1 + q\theta_2) \biggr)\alpha_1^i \alpha_2^j \alpha_3^k \alpha_4^m \\
	y(t,\alpha_1,\alpha_2,\alpha_3,\alpha_4) = \sum e^{(i-j)\theta_3} \biggr(y_{ijkm}^{pq} \cos(p\theta_1 + q\theta_2) + \bar{y}_{ijkm}^{pq} \sin(p\theta_1 + q\theta_2) \biggr)\alpha_1^i \alpha_2^j \alpha_3^k \alpha_4^m \\
	z(t,\alpha_1,\alpha_2,\alpha_3,\alpha_4) = \sum e^{(i-j)\theta_3} \biggr(z_{ijkm}^{pq} \cos(p\theta_1 + q\theta_2) + \bar{z}_{ijkm}^{pq} \sin(p\theta_1 + q\theta_2) \biggr)\alpha_1^i \alpha_2^j \alpha_3^k \alpha_4^m
\end{aligned}	
\label{eqn:series}
\end{equation}
In these expressions, $\theta_1 = \omega t + \phi_1$, $\theta_2 = \nu t + \phi_2$, and $\theta_3 = \lambda t$ where $\omega$ is the in-plane frequency, $\nu$ is the out-of-plane frequency, $\lambda$ encompasses the hyperbolic nature of the system corresponding to the manifolds, and $\phi_1$ and $\phi_2$ are arbitrary phases. The frequencies are also solved for as a sum of coefficients.
\begin{equation}
\begin{aligned}
	\omega = \sum \omega_{ijkm} \alpha_1^i \alpha_2^j \alpha_3^k \alpha_4^m \\
	\nu = \sum \nu_{ijkm} \alpha_1^i \alpha_2^j \alpha_3^k \alpha_4^m \\ 
	\lambda = \sum \lambda_{ijkm} \alpha_1^i \alpha_2^j \alpha_3^k \alpha_4^m
\end{aligned}	
\end{equation}
The Lindstedt-Poincar\'e method solves for the coefficients of the solutions iteratively up to a user-defined order $n$. The sum of $i,j,k,m$ defines the order $n$ of the solution and $p \leq k$ (mod 2) and $|q| \leq m$ (mod 2). Additionally, when $p=0$, $q \leq m$. The amplitudes $\alpha_1$, $\alpha_2$, $\alpha_3$, and $\alpha_4$ parameterize a spacecraft's motion about an AEP in the following manner.
\begin{enumerate}[noitemsep,topsep=0pt]
	\item When $\alpha_1 \neq 0$ and $\alpha_2 = 0$, Eq.~\eqref{eqn:series} yields unstable manifolds.
	\item When $\alpha_1 = 0$ and $\alpha_2 \neq 0$, Eq.~\eqref{eqn:series} yields stable manifolds.
	\item When $\alpha_1,\alpha_2 \neq 0$ and sign$(\alpha_1) \neq $ sign$(\alpha_2)$, Eq.~\eqref{eqn:series} yields transit trajectories that move the spacecraft from one side of the AEP to the other.
	\item When $\alpha_1,\alpha_2 \neq 0$ and sign$(\alpha_1) = $ sign$(\alpha_2)$, Eq.~\eqref{eqn:series} yields non-transit trajectories that move the spacecraft from one side of the AEP and returns it to the same side.
	\item When $\alpha_1,\alpha_2= 0 $ and either $\alpha_3$ or $\alpha_4$ is nonzero, Eq.~\eqref{eqn:series} produces a periodic orbit.
\end{enumerate}

The Lindstedt-Poincar\'e method is an iterative and recursive technique, employing the results from the previous order to solve the current order. This requires that the process be initialized with the zeroth-order and first-order solutions from the linearized equations of motion. The zeroth-order and first-order solutions of the coordinate series are
\begin{equation}
\begin{aligned}
	&x_{1000}^{00} = 1, \quad x_{0100}^{00} = 1, \quad x_{0010}^{10} = 1, \quad \bar{x}_{0010}^{10} = k_{14} \\
	&x_{0001}^{01} = k_5, \quad \bar{x}_{0001}^{01} = k_6, \quad x_{0000}^{00} = k_{15} \\
	&y_{1000}^{00} = k_1, \quad y_{0100}^{00} = k_2, \quad y_{0010}^{10} = k_7, \quad \bar{y}_{0010}^{10} = k_{8} \\
	&y_{0001}^{01} = k_9, \quad \bar{y}_{0001}^{01} = k_{10}, \quad y_{0000}^{00} = k_{16}  \\
	&z_{1000}^{00} = k_3, \quad z_{0100}^{00} = k_4, \quad z_{0010}^{10} = k_{11}, \quad \bar{z}_{0010}^{10} = k_{12} \\
	&z_{0001}^{01} = 1, \quad \bar{z}_{0001}^{01} = k_{13}, \quad z_{0000}^{00} = k_{17} . 
\end{aligned}	
\end{equation}
The zeroth-order solutions of the frequency series are
\begin{equation}
\omega_{0000} = \omega_0, \quad \nu_{0000} = \nu_0, \quad \lambda_{0000} = \lambda_0.
\end{equation}
The order of the solution is defined as $n = i + j + k + m$. At each order $n$, the known series solutions up to order $n-1$ are substituted into the power series of the equations of motion and the terms corresponding to order $n$ are held. These known values are denoted by $r_{ijkm}^{pq}$, $\bar{r}_{ijkm}^{pq}$, $s_{ijkm}^{pq}$, $\bar{s}_{ijkm}^{pq}$, $t_{ijkm}^{pq}$, and $\bar{t}_{ijkm}^{pq}$ --- this essentially groups all like terms together for each coordinate series.

The unknown coefficients of order $n$ are solved via the system of linear equations:
\begin{equation}
M \mathbf{X} + \boldsymbol{\delta} = \mathbf{b} + \mathbf{c}
\end{equation}
where 
\begin{equation}
	M = \left [ 
	\begin{array}{cccccc}
		\xi- \Omega_{xx}^* & 2\zeta\Psi & -2\zeta-\Omega_{xy}^* & -2\Psi & -\Omega_{xz}^* & 0 \\
		-2\zeta\Psi & \xi - \Omega_{xx}^* & 2\Psi & -2\zeta-\Omega_{xy}^* & 0 &  -\Omega_{xz}^* \\
		2\zeta - \Omega_{yx}^* & 2\Psi & \xi - \Omega_{yy}^* & 2\zeta\Psi & -\Omega_{yz}^* & 0 \\ 
		-2\Psi & 2\zeta - \Omega_{yx}^* & -2\zeta\Psi & \xi - \Omega_{yy}^* & 0 & -\Omega_{yz}^* \\
		-\Omega_{zx}^* & 0 & -\Omega_{zy}^* & 0 & \xi - \Omega_{zz}^* & 2\zeta\Psi \\
		0 & -\Omega_{zx}^* & 0 & -\Omega_{zy}^* & -2\zeta\Psi & \xi - \Omega_{zz}^*
	\end{array}
	\right ]
\end{equation}
$\Psi = p\omega_0 + q\nu_0$, $\zeta = (i-j)\lambda_0$, and $\xi = \zeta^2  - \Psi^2$, 
\begin{equation}
	\mathbf{X} = \left[
	\begin{array}{c}
		x_{ijkm}^{pq} \\
		\bar{x}_{ijkm}^{pq} \\
		y_{ijkm}^{pq} \\
		\bar{y}_{ijkm}^{pq} \\
		z_{ijkm}^{pq} \\
		\bar{z}_{ijkm}^{pq} \\
	\end{array}
	\right], \quad \text{and } \quad
	\mathbf{b} = \left[
	\begin{array}{c}
		r_{ijkm}^{pq} \\
		\bar{r}_{ijkm}^{pq} \\
		s_{ijkm}^{pq} \\
		\bar{s}_{ijkm}^{pq} \\
		t_{ijkm}^{pq} \\
		\bar{t}_{ijkm}^{pq} \\
	\end{array}
	\right].
\end{equation}
The contents of $\boldsymbol{\delta}$ and $\mathbf{c}$ vary with the set of indices ($i,\, j, \, k, \, m,\, p, \, q$) as follows.
\begin{enumerate}[noitemsep,topsep=0pt]
	\item When $p = q = 0$
	\begin{enumerate}
		\item and $i - 1 = j$, $\boldsymbol{\delta} = [2(\lambda_0 - k_1), 0, 2(k_1\lambda_0+1), 0, 2 k_3 \lambda_0]^T \lambda_{i-1jkm}$ \newline and $\mathbf{c} = -\bar{\Lambda}_{i-1jkm}[1, 0, k_1, 0, k_3, 0]^T$.
		\item and $i = j-1$,  $\boldsymbol{\delta} = [2(\lambda_0 + k_2), 0, 2(k_2\lambda_0-1), 0, 2 k_4 \lambda_0]^T \lambda_{ij-1km}$ \newline and $\mathbf{c} = -\bar{\Lambda}_{ij-1km}[1, 0, k_2, 0, k_4, 0]^T$.
	\end{enumerate}
	\item When $p = 1$, $q = 0$, and $i=j$, $\boldsymbol{\delta}=[-2(k_8 + \omega_0), 2(k_7-k_{14}\omega_0), 2(k_{14} -k_7\omega_0), -2(k_8\omega_0+1), -2k_{11}\omega_0, -2k_{12}\omega_0]^T\omega_{iik-1m}$ and $\mathbf{c} = \bar{\Omega}_{iik-1m} [1, k_{14}, k_7, k_8, k_{11}, k_{12} ]^T $.
	\item When $p = 0$, $q = 1$, and $i=j$, $\boldsymbol{\delta}=[-2(k_5\nu_0 + k_{10}), -2(k_6\nu_0 - k_9), -2(k_9\nu_0-k_6), -2(k_{10}\nu_0+k_5), -2\nu_0, -2k_{13}\nu_0]^T\nu_{iikm-1}$ and $\mathbf{c} = \bar{N}_{iikm-1} [k_5, k_6, k_9, k_{10}, 1,k_{13} ]^T $.
	\item Otherwise, $\boldsymbol{\delta} = \mathbf{0}$ and $\mathbf{c} = \mathbf{0}$.
\end{enumerate}
In cases (1)-(3), the unknown frequency is appended to $\mathbf{X}$ and $\boldsymbol{\delta}$ is concatenated with $M$. Therefore, a pseudoinverse of the new $M$ matrix is required to yield the unknown coefficients. This is not required in case (4) since $M$ is invertible.

The resulting coordinate series from the Lindstedt-Poincar\'e method can be differentiated with respect to time to produce the velocity series. In this manner, a full-state approximation for CR3BP with SRP can be obtained.

\section{Results}
This section considers results yielded from periodic orbits and manifolds about the AEP analogous with and in the vicinity of Sun-Earth $L_2$. For all cases considered herein, $\beta = 0.002$ is the expected sail lightness number associated with the HabEx mission. The non-dimensional mass parameter of the Sun-Earth system used is $\mu = 3.0026053634189284 \times 10^{-6}$. The periodic orbits that exist for non-Hamiltonian SRP forces break from the characterizations associated with classic CR3BP orbits. Namely, the resulting periodic orbits do not generally share the same properties as Lyapunov, vertical, Lissajous, or Halo orbits; a comparison to these types of orbits will be discussed for the examples considered.

Three sets of attitude angles are considered here: ($\alpha$, $\gamma$) = ($80^\circ$, $0^\circ$), ($0^\circ$, $40^\circ$), and ($80^\circ$, $40^\circ$). Using the Lindstedt-Poincar\'e method discussed above, seventh-order approximate solutions were constructed and analyzed for the three cases. The accuracy of the resulting periodic orbits is found by comparing the half orbit produced against that yielded from the original equations of motion in Eq.~\eqref{eqn:eoms}\footnotemark{}\footnotetext{The half orbit is used for an accuracy comparison rather than a full orbit due to the lack of a general differential corrector scheme that can refine the initial condition produced by the approximate solution, which is fed into an ODE solver for the original equations of motion. When considering the CR3BP analogs, it is also seen that, without a differential corrector, the two orbits start to diverge significantly after the half-period mark.}. For each set, stable and unstable manifolds are found, sample periodic orbits (analogous to Lyapunov, vertical, and Lissajous) are displayed, and the results of a numerical study to find the position error with respect to truth are presented.

\subsection{Approximate Solutions for $\alpha = 80^\circ$, $\gamma = 0^\circ$}
\begin{figure}[!htp]
\centering
\subfloat[Lyapunov-like]{%
  \includegraphics[width=45mm]{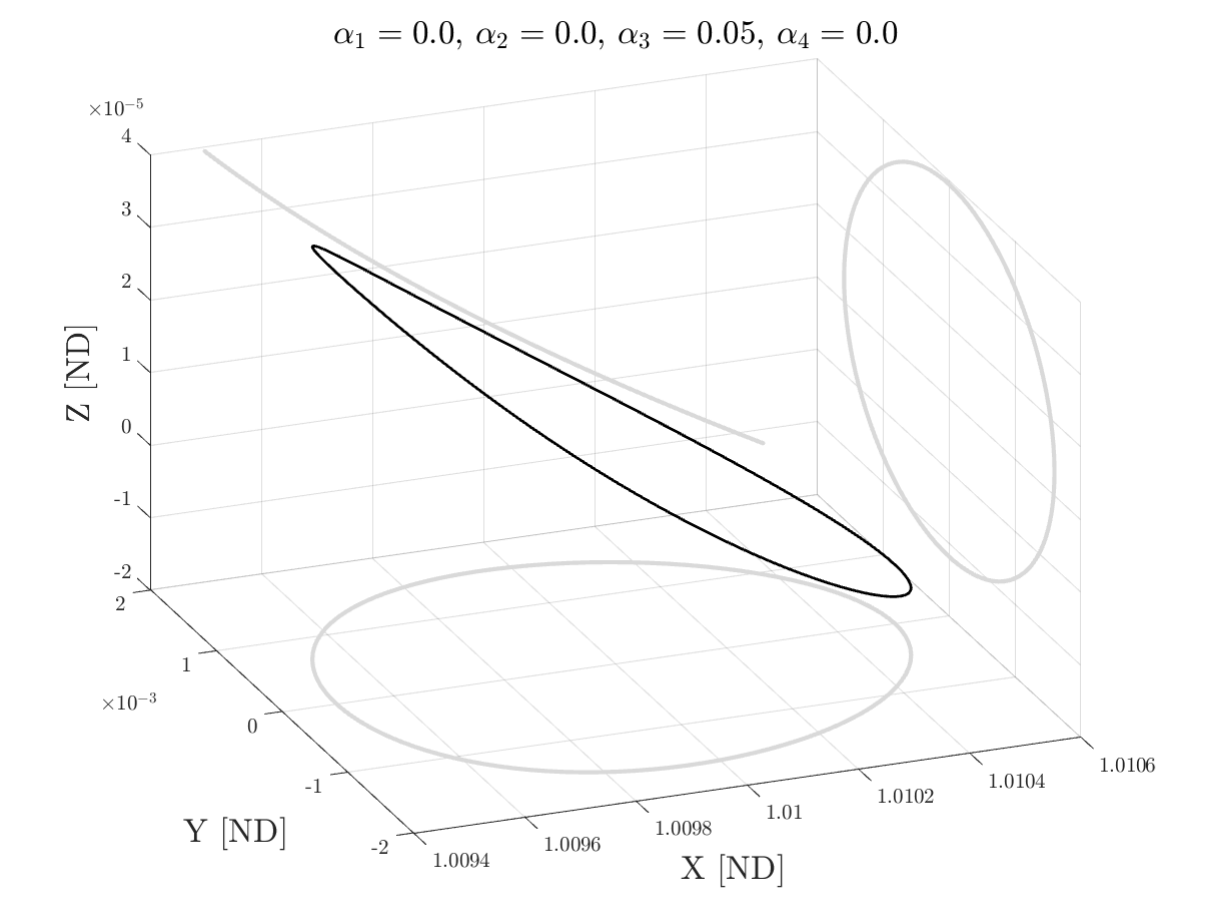}\label{fig:alpha80Lyapunov}%
}\quad
\subfloat[Vertical-like]{%
  \includegraphics[width=45mm]{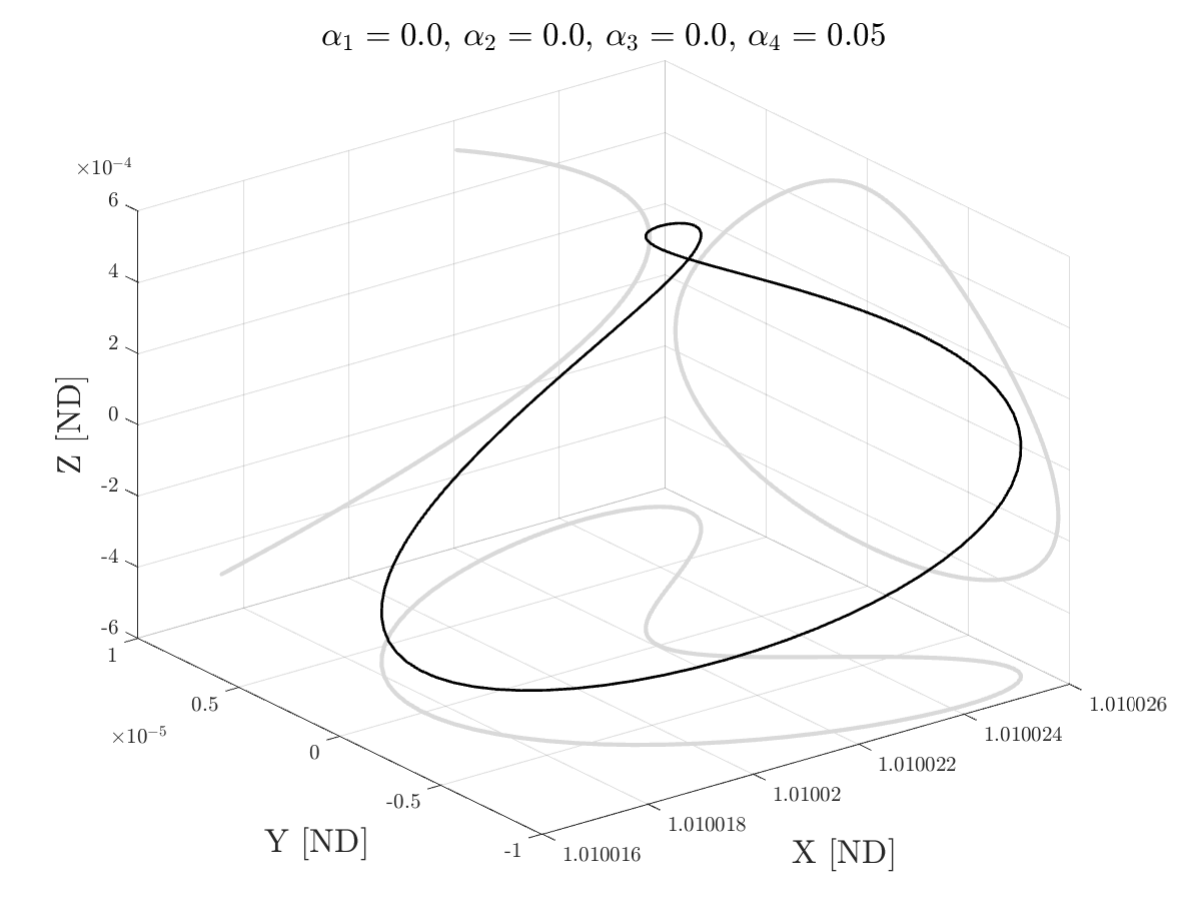}\label{fig:alpha80Vertical}
}\quad
\subfloat[Lissajous-like]{%
  \includegraphics[width=45mm]{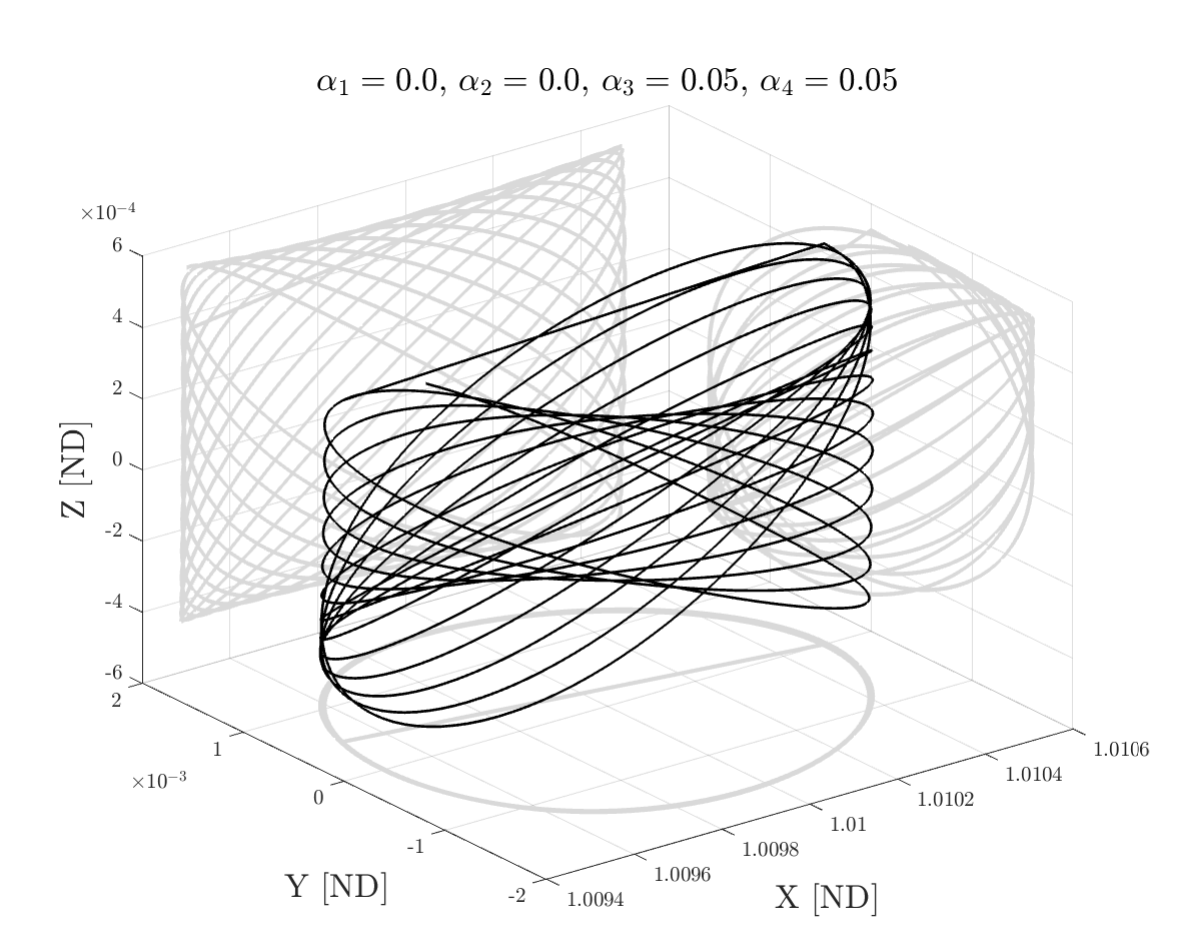}\label{fig:alpha80Lissajous}
}
\caption{Periodic Orbits for $\alpha = 80^\circ$, $\gamma = 0^\circ$}
\label{fig:alpha80}
\end{figure}

Figure \ref{fig:alpha80} shows examples of the three different types of periodic orbits for ($\alpha$, $\gamma$) = ($80^\circ$, $0^\circ$), produced using the seventh-order approximation. 
 \newline
For this set of attitude angles, the position of $H_2$ is
\begin{equation*}
	\left[
	\begin{array}{ccc}
		 1.0100319725242741, & 0.0, & 1.4769123813475747\times 10^{-5}
	\end{array}
	\right ]
\end{equation*}
in non-dimensional units. The nonzero $z$-component of the position indicates that the AEP is not collinear with the $x$-axis in the rotating frame. Therefore, it is impossible to produce purely planar, periodic orbits (i.e. Lyapunov orbits), which, for CR3BP, are found using a nonzero $\alpha_3$ value and $\alpha_4 = 0$. With this set of attitude angles, $\alpha_3 \neq 0$ and $\alpha_4 = 0$ produces a three-dimensional periodic orbit, shown in Figure \ref{fig:alpha80Lyapunov}. Similarly, CR3BP with SRP does not produce the traditional folded-over figure-eight shape for vertical orbits (which occurs when $\alpha_3 = 0$ and $\alpha_4 \neq 0$). However, for this set of attitude angles, these amplitude values still produce a folded-over orbit, seen in Figure \ref{fig:alpha80Vertical}. Lastly, when both $\alpha_3$ and $\alpha_4$ are nonzero, the characteristic Lissajous shape is retained.

\begin{figure}[!htp]
\centering
\subfloat[Stable manifolds]{%
  \includegraphics[width=70mm]{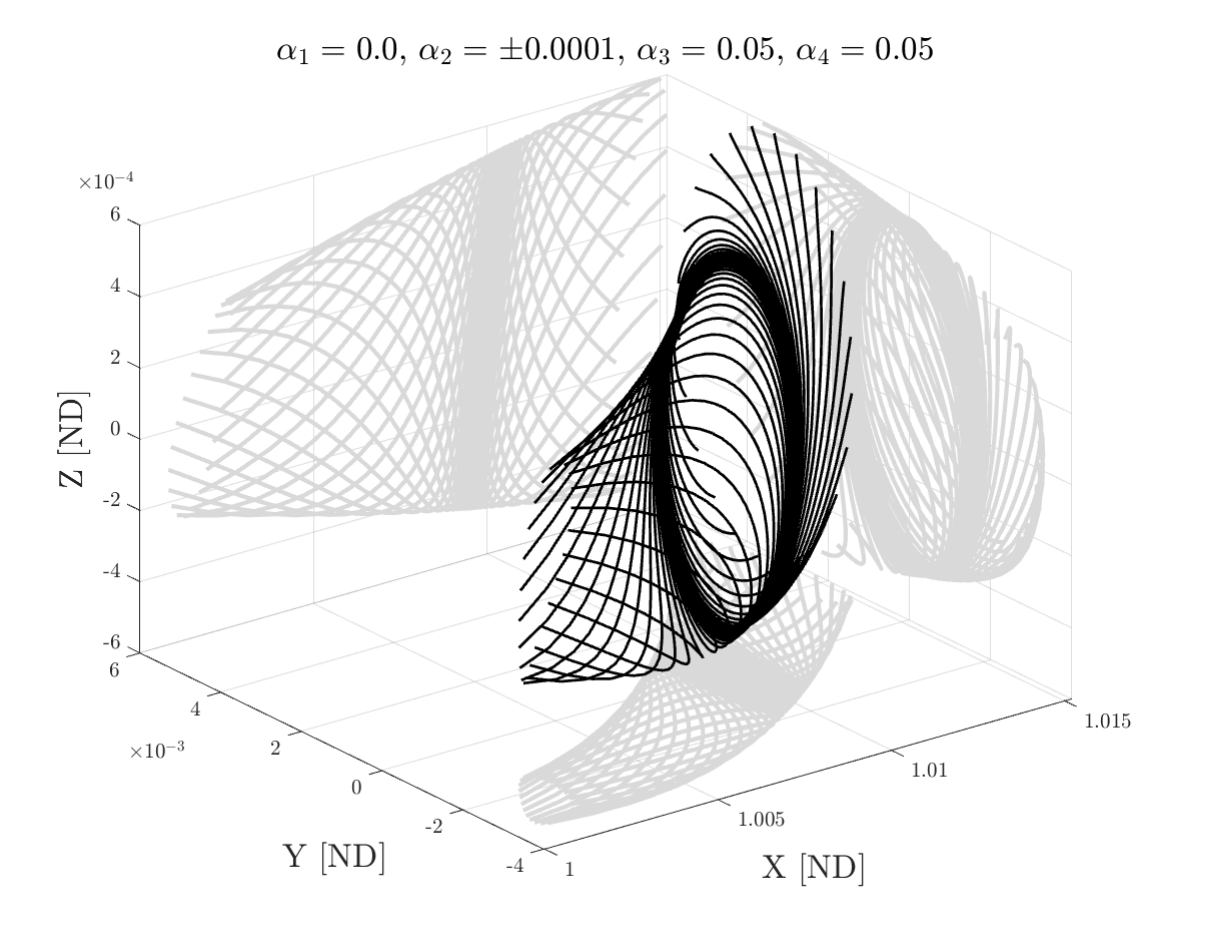}%
}\quad
\subfloat[Unstable manifolds]{%
  \includegraphics[width=70 mm]{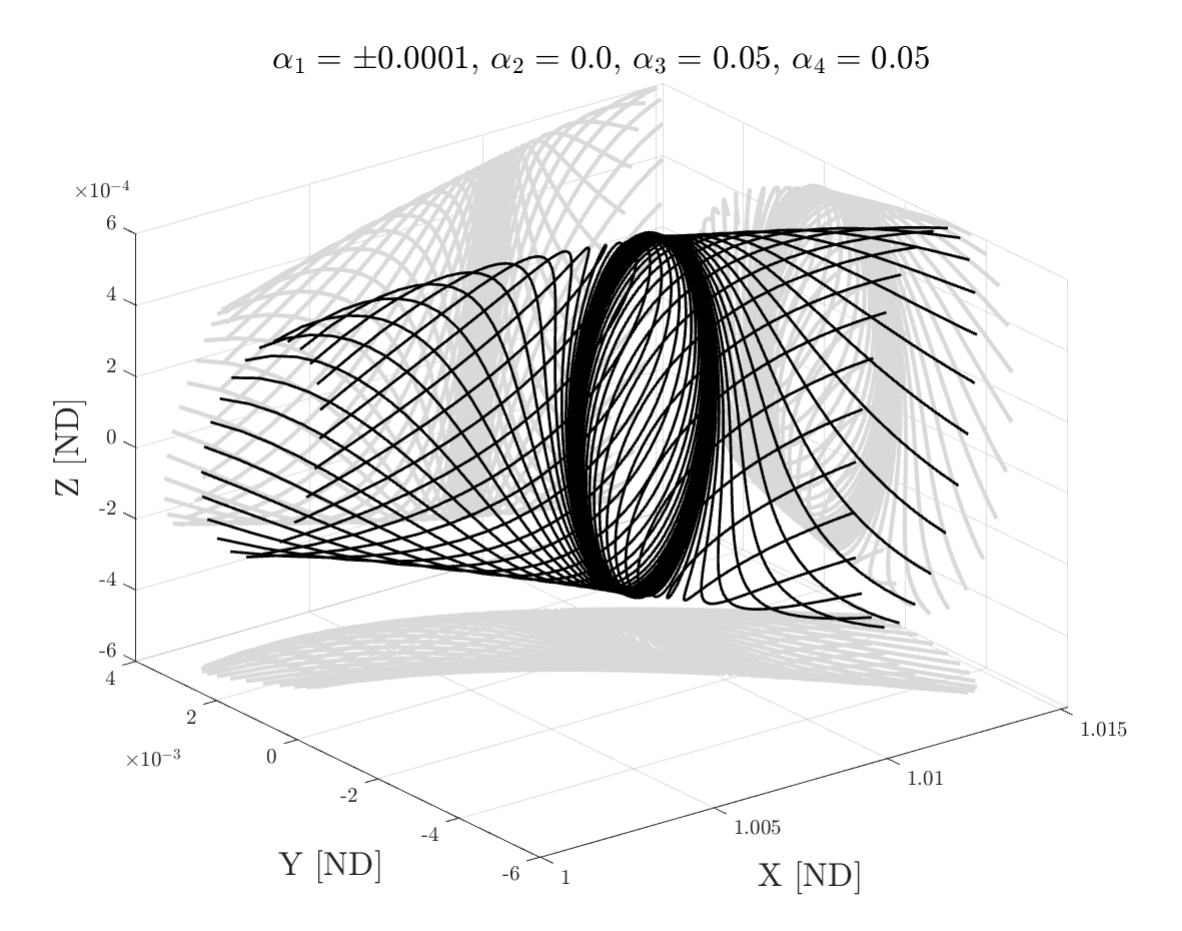}
}
\caption{Manifolds for $\alpha = 80^\circ$, $\gamma = 0^\circ$}
\label{fig:results1}
\end{figure}

In addition to periodic orbits, the approximate solutions may be used to produce stable and unstable manifolds, which lead to and from the AEP, respectively. In Figure \ref{fig:results1}, $\alpha_3 = \alpha_4 = 0.05             $ for both plots. In the lefthand plot, $\alpha_1 = 0 $ and $\alpha_2 =  \pm 1 \times 10^{-4}$. In the righthand plot, $\alpha_1 = \pm 1 \times 10^{-4}$ and $\alpha_2 = 0$. Similar to the periodic orbits, these manifolds uniquely exist for \textit{this AEP}, suggesting that a change in one or both of the attitude angles can produce a unique set of manifolds that may better suit a specific mission's purpose.

\subsection{Approximate Solutions for $\alpha = 0^\circ$, $\gamma = 40^\circ$}
\begin{figure}[!htp]
\centering
\subfloat[Lyapunov-like]{%
  \includegraphics[width=45mm]{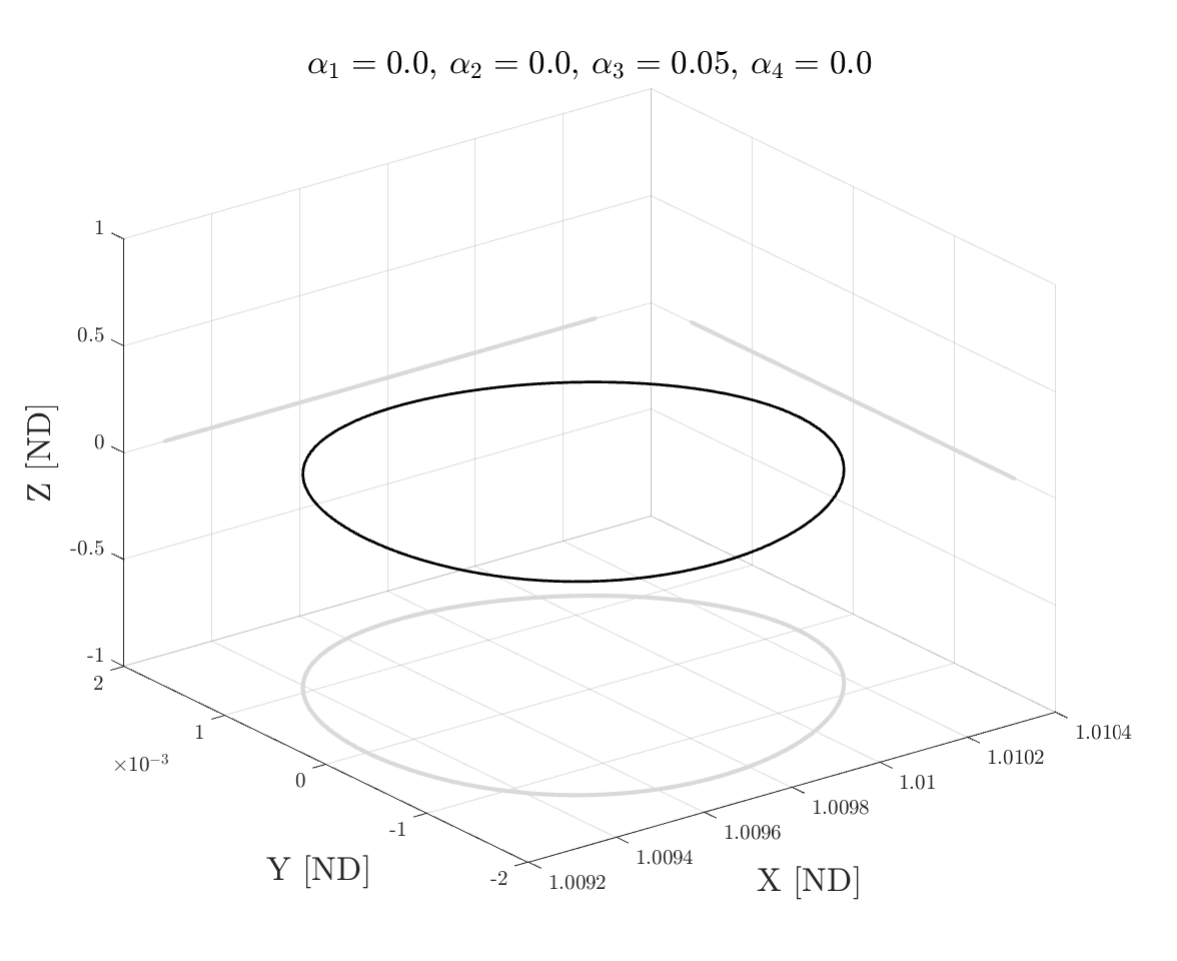}\label{fig:gamma40Lyapunov}%
}\quad
\subfloat[Vertical-like]{%
  \includegraphics[width=45mm]{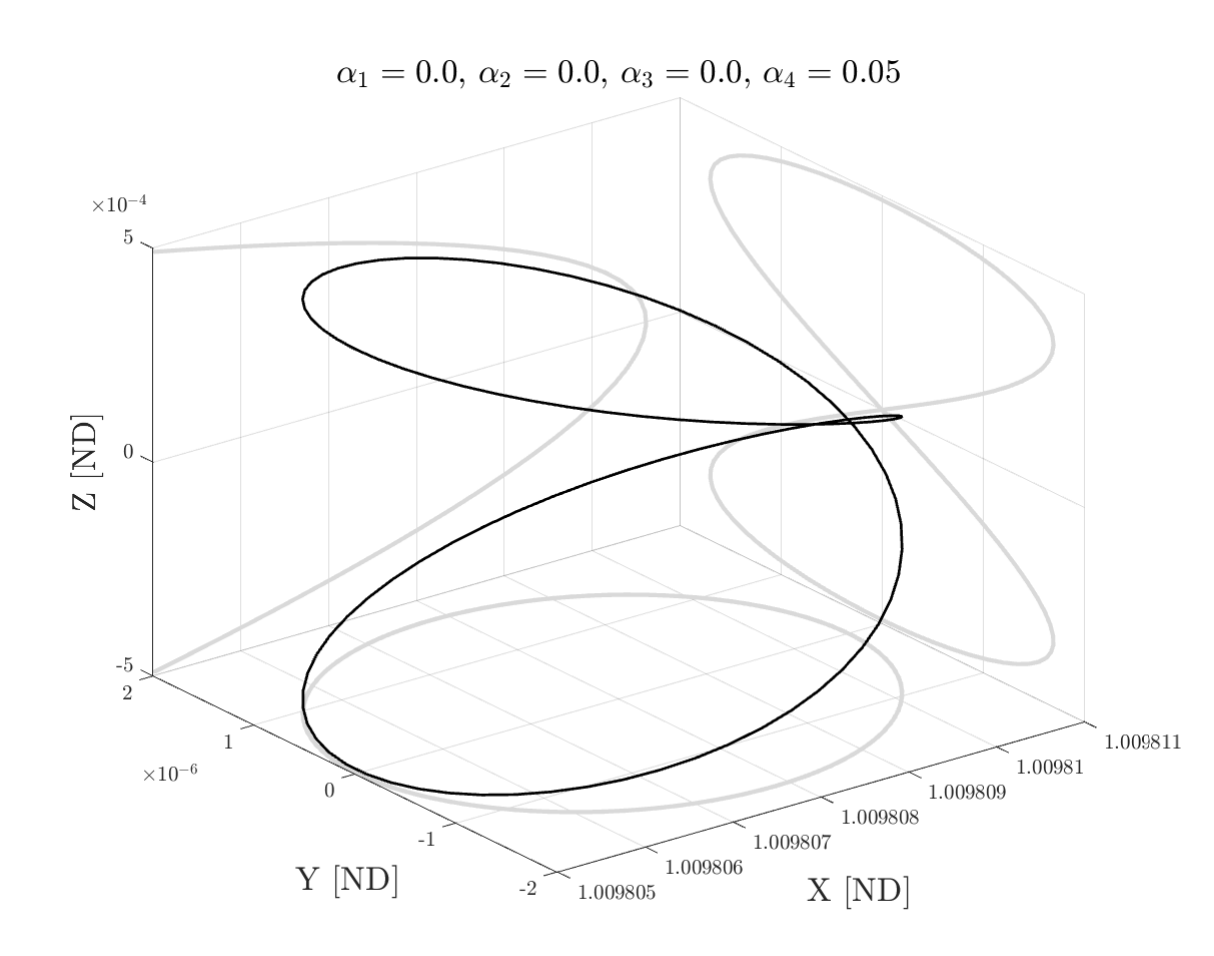}\label{fig:gamma40Vertical}
}\quad
\subfloat[Lissajous-like]{%
  \includegraphics[width=45mm]{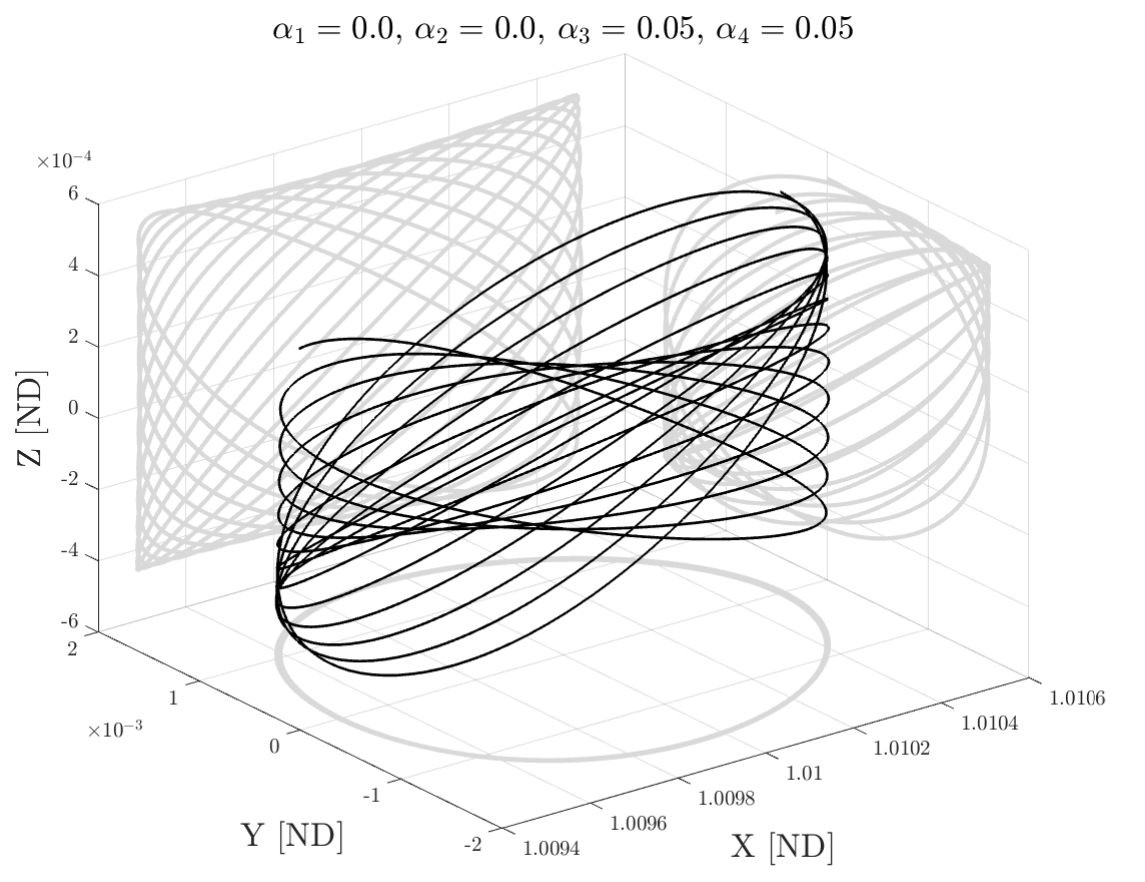}\label{fig:gamma40Lissajous}
}
\caption{Periodic Orbits for $\alpha = 0^\circ$, $\gamma = 40^\circ$}
\label{fig:gam40}
\end{figure}

Figure \ref{fig:gam40} shows examples of the the three different types of periodic orbits for ($\alpha$, $\gamma$) = ($0^\circ$, $40^\circ$), produced using the seventh-order approximation. 
\newline
For this set of attitude angles, the position of $H_2$ is
\begin{equation*}
	\left[
	\begin{array}{ccc}
		 1.009817129039308, & 0.0, & 0.0
	\end{array}
	\right ]
\end{equation*}
in non-dimensional units. Since the AEP is collinear with the $x$-axis in the rotating frame, the characteristic shapes of the orbits from CR3BP are maintained. Physically, this is intuitive when it is recognized that $\alpha = 0^\circ$ means the starshade is exactly perpendicular to the Sun-line, acting as a conservative force now since it is path-independent.

\begin{figure}[!htp]
\centering
\subfloat[Stable manifolds]{%
  \includegraphics[width=70mm]{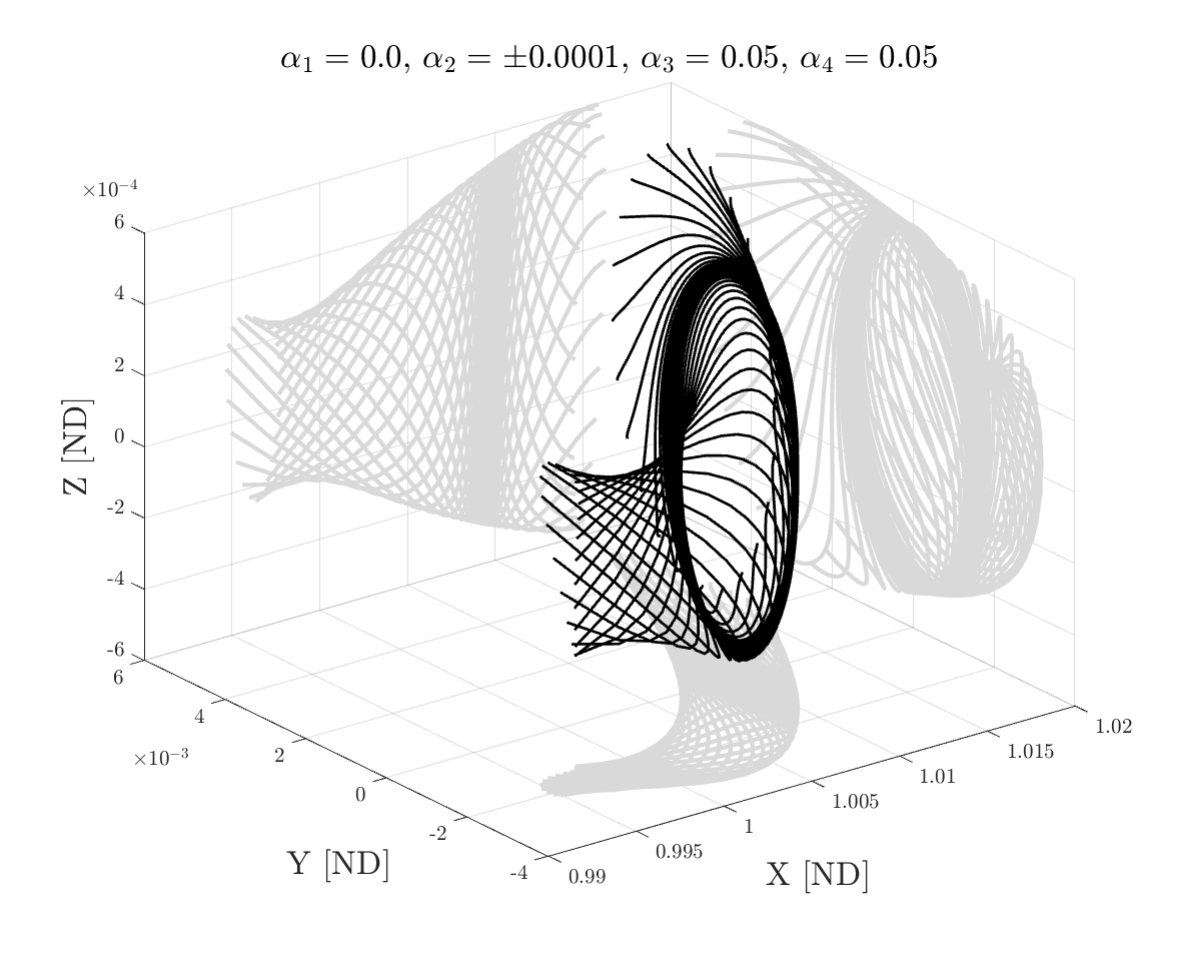}%
}\quad
\subfloat[Unstable manifolds]{%
  \includegraphics[width=70 mm]{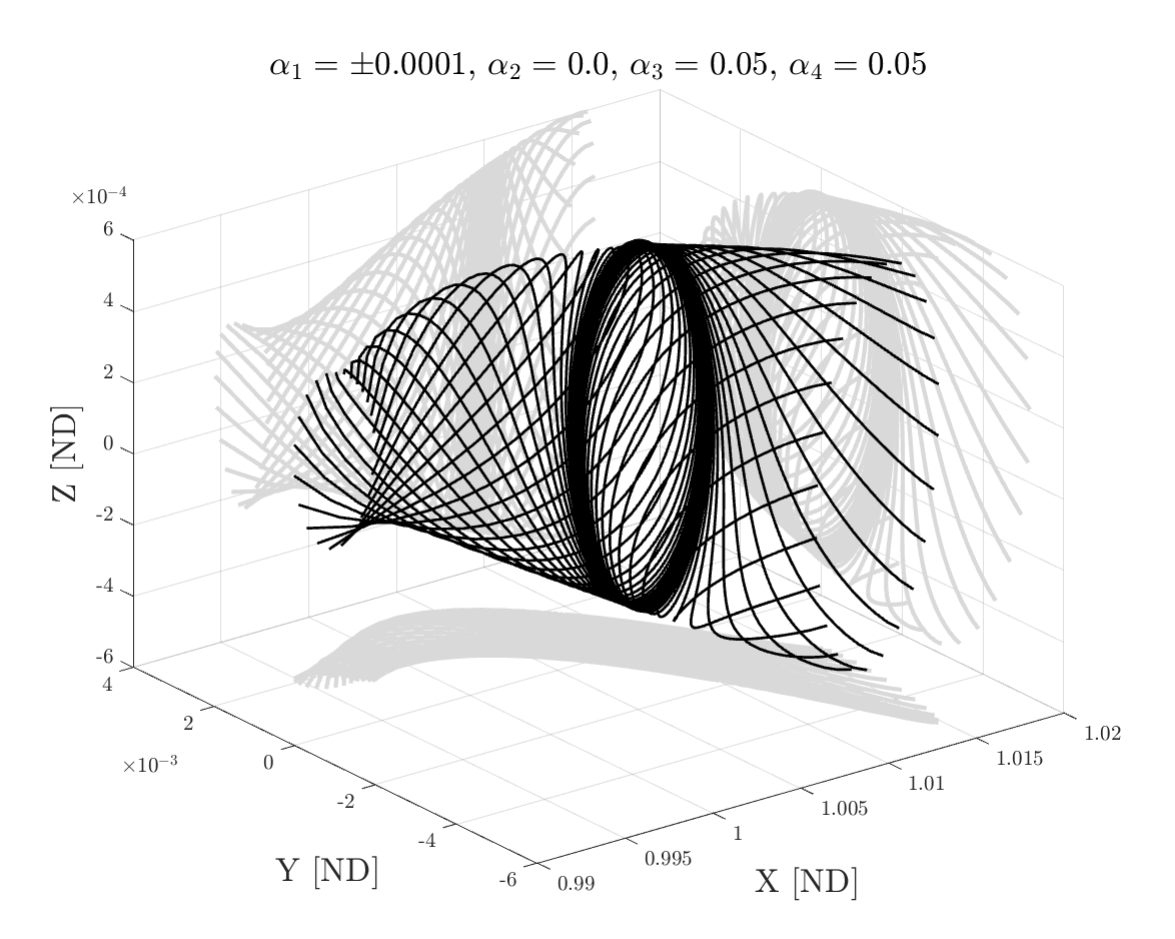}
}
\caption{Manifolds for $\alpha = 0^\circ$, $\gamma = 40^\circ$}
\label{fig:results2}
\end{figure}

In addition to periodic orbits, the approximate solutions may be used to produce stable and unstable manifolds, which lead to and from the AEP, respectively. In Figure \ref{fig:results2}, $\alpha_3 = \alpha_4 = 0.05$ for both plots. In the lefthand plot, $\alpha_1 = 0 $ and $\alpha_2 =  \pm 1 \times 10^{-4}$. In the righthand plot, $\alpha_1 = \pm 1 \times 10^{-4}$ and $\alpha_2 = 0$. Similar to the periodic orbits, these manifolds uniquely exist for \textit{this AEP}, suggesting that a change in one or both of the attitude angles can produce a unique set of manifolds that may better suit a specific mission's purpose.

\subsection{Approximate Solutions for $\alpha = 80^\circ$, $\gamma = 40^\circ$}
\begin{figure}[!htp]
\centering
\subfloat[Lyapunov-like]{%
  \includegraphics[width=45mm]{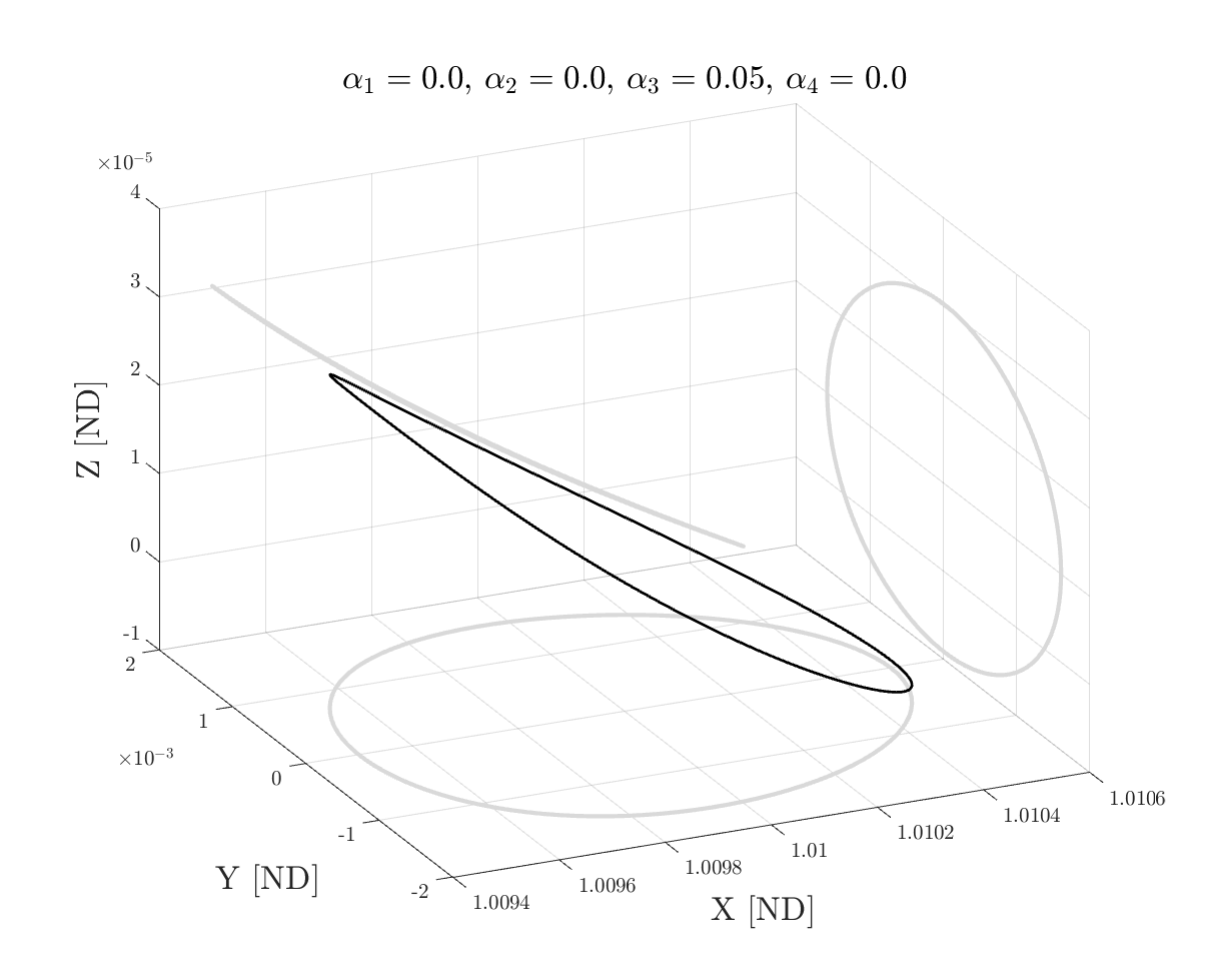}\label{fig:alpha80gamma40Lyapunov}%
}\quad
\subfloat[Vertical-like]{%
  \includegraphics[width=45mm]{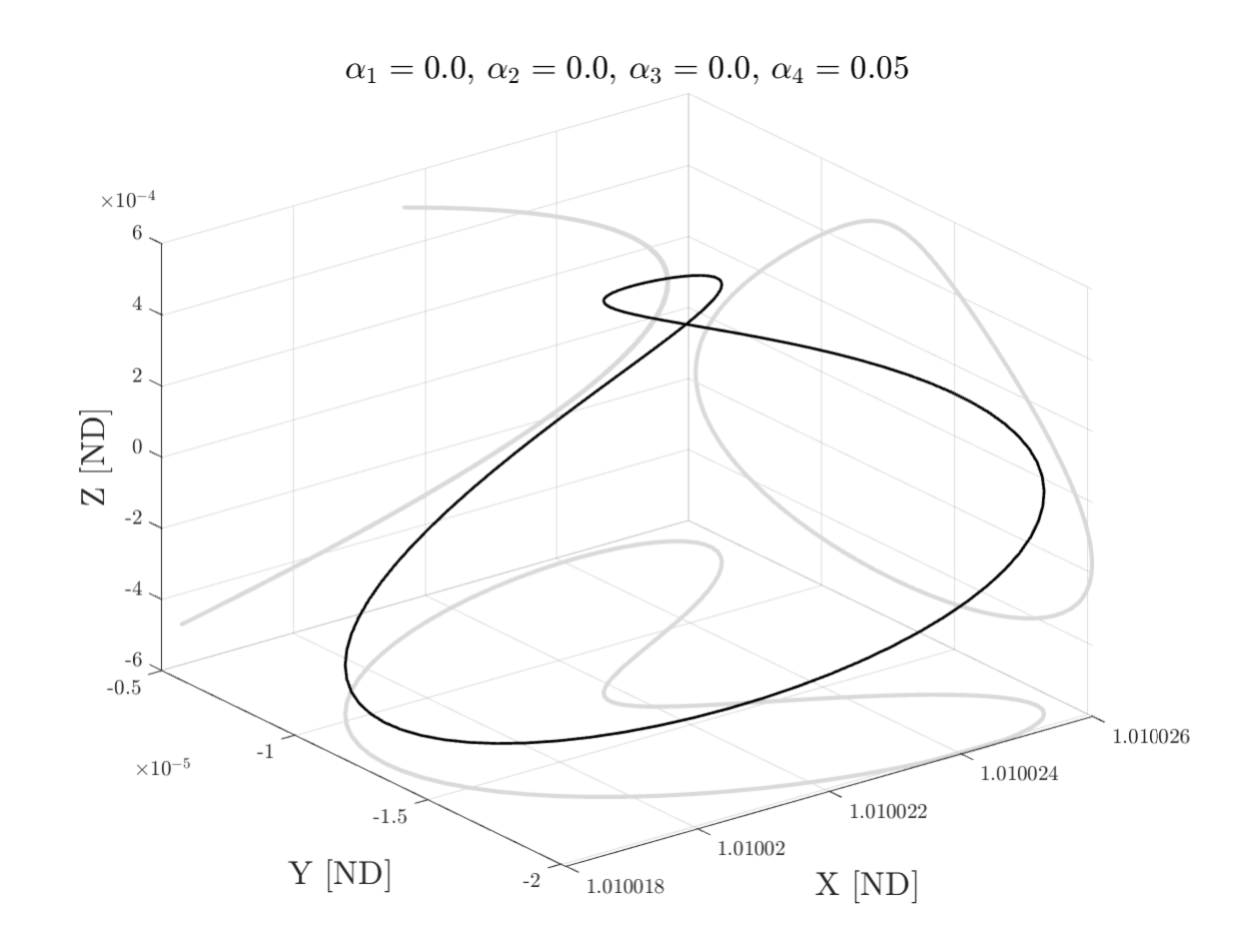}\label{fig:alpha80gamma40Vertical}
}\quad
\subfloat[Lissajous-like]{%
  \includegraphics[width=45mm]{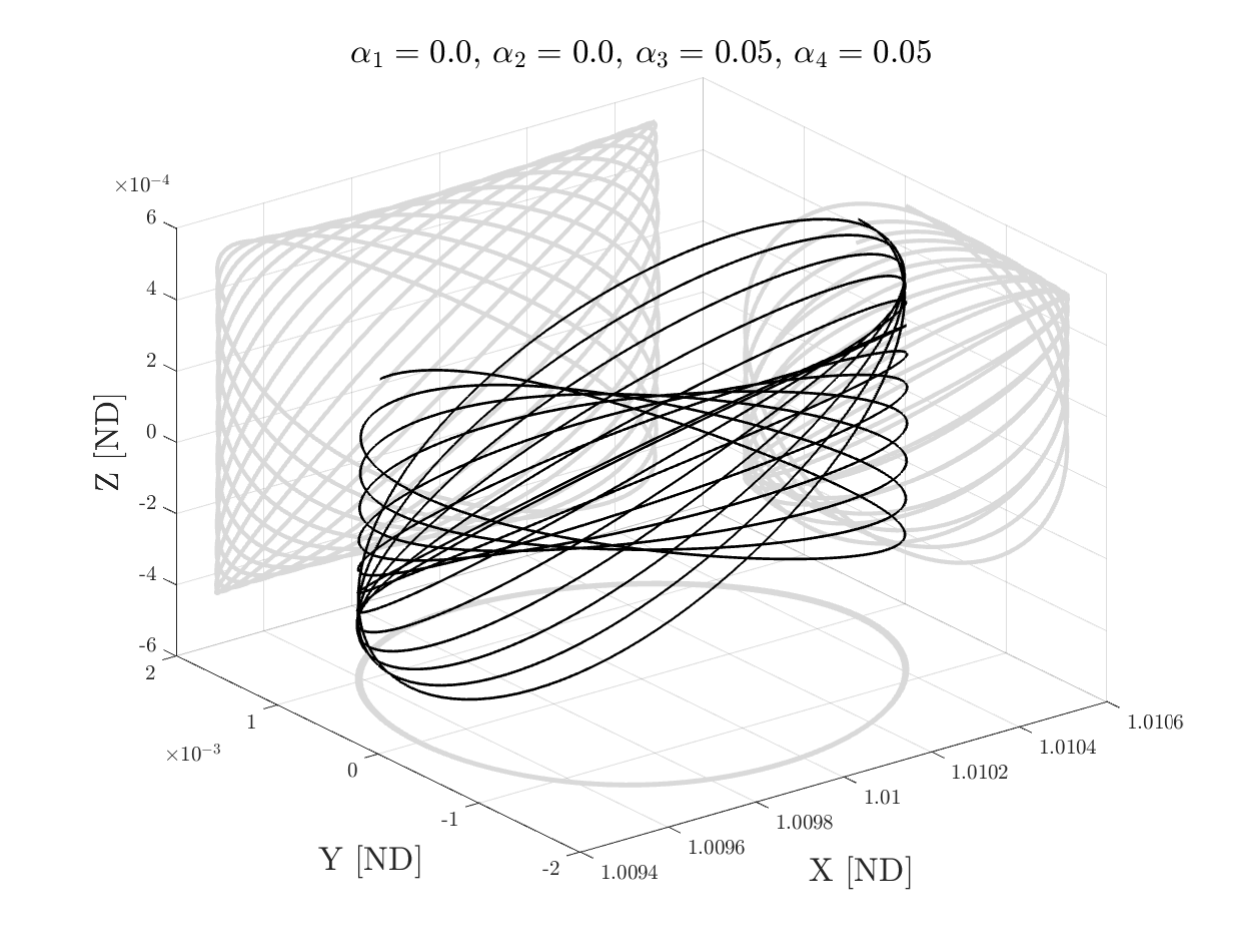}\label{fig:alpha80gamma40Lissajous}
}
\caption{Periodic Orbits for $\alpha = 80^\circ$, $\gamma = 40^\circ$}
\label{fig:alp80gam40}
\end{figure}

Figure \ref{fig:alp80gam40} shows examples of the three different types of periodic orbits for ($\alpha$, $\gamma$) = ($80^\circ$, $40^\circ$), produced using the seventh-order approximation. 
\newline
For this set of attitude angles, the position of $H_2$ is
\begin{equation*}
	\left[
	\begin{array}{ccc}
		 1.0100319689420738, & -1.2720500232390416\times 10^{-5}, & 1.1313805251204233\times 10^{-5}
	\end{array}
	\right ]
\end{equation*}
in non-dimensional units. The position's nonzero $y$-component and $z$-component indicate that the AEP is not collinear with the $x$-axis in the rotating frame. This indicates that there do not exist purely planar, periodic orbits like in CR3BP, which are found using a nonzero $\alpha_3$ value and $\alpha_4 = 0$. For this amplitude configuration and set of attitude angles, a three-dimensional periodic orbit is produced, shown in Figure \ref{fig:alpha80gamma40Lyapunov}. Similar to the ($\alpha$, $\gamma$) = ($80^\circ$, $0^\circ$) case, the classic vertical orbits are not possible, but there still exist orbits distantly analogous to its shape (Figure \ref{fig:alpha80gamma40Vertical}). Lastly, when both $\alpha_3$ and $\alpha_4$ are nonzero, the characteristic Lissajous shape is retained, as shown in Figure \ref{fig:alpha80gamma40Lissajous}. 
\begin{figure}[!htp]
\centering
\subfloat[Stable manifolds]{%
  \includegraphics[width=70mm]{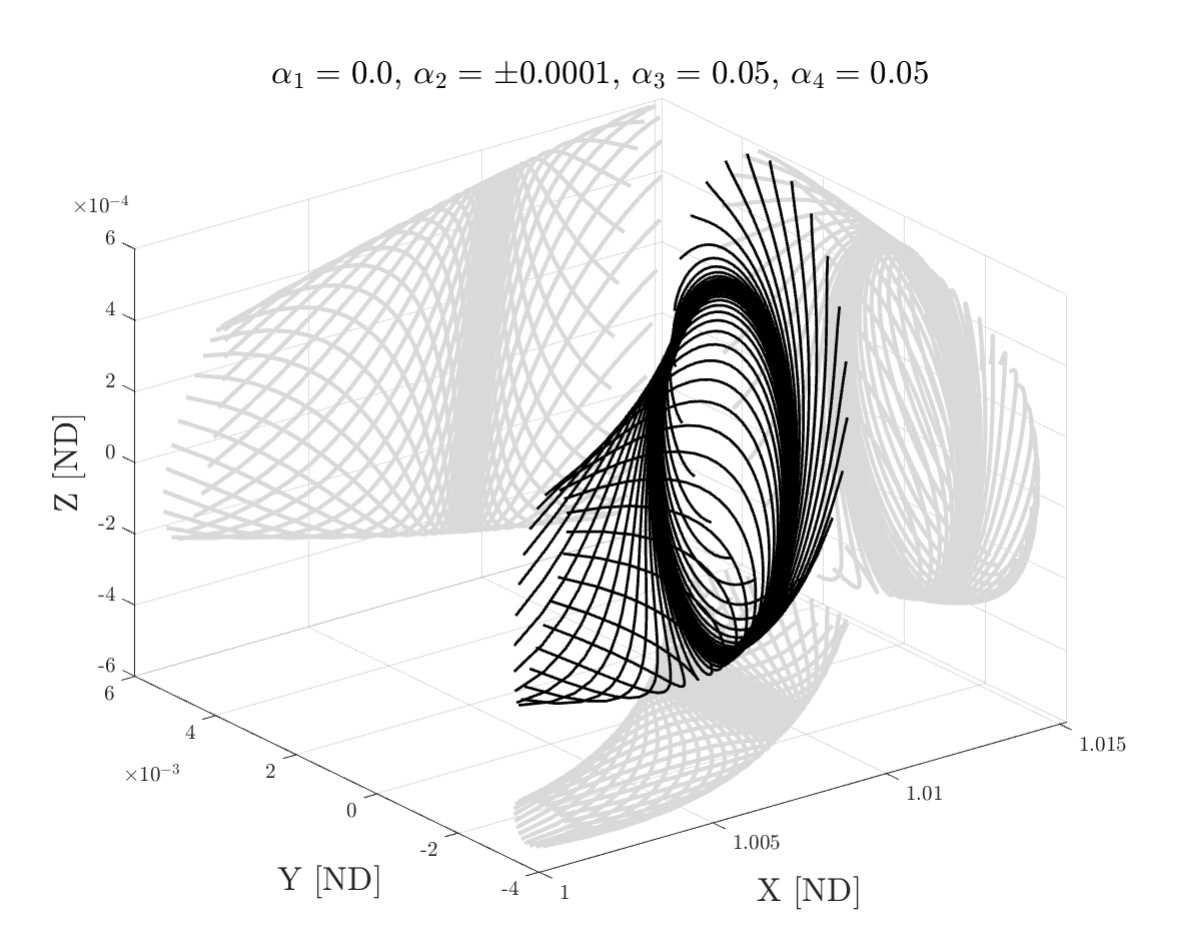}%
}\quad
\subfloat[Unstable manifolds]{%
  \includegraphics[width=70 mm]{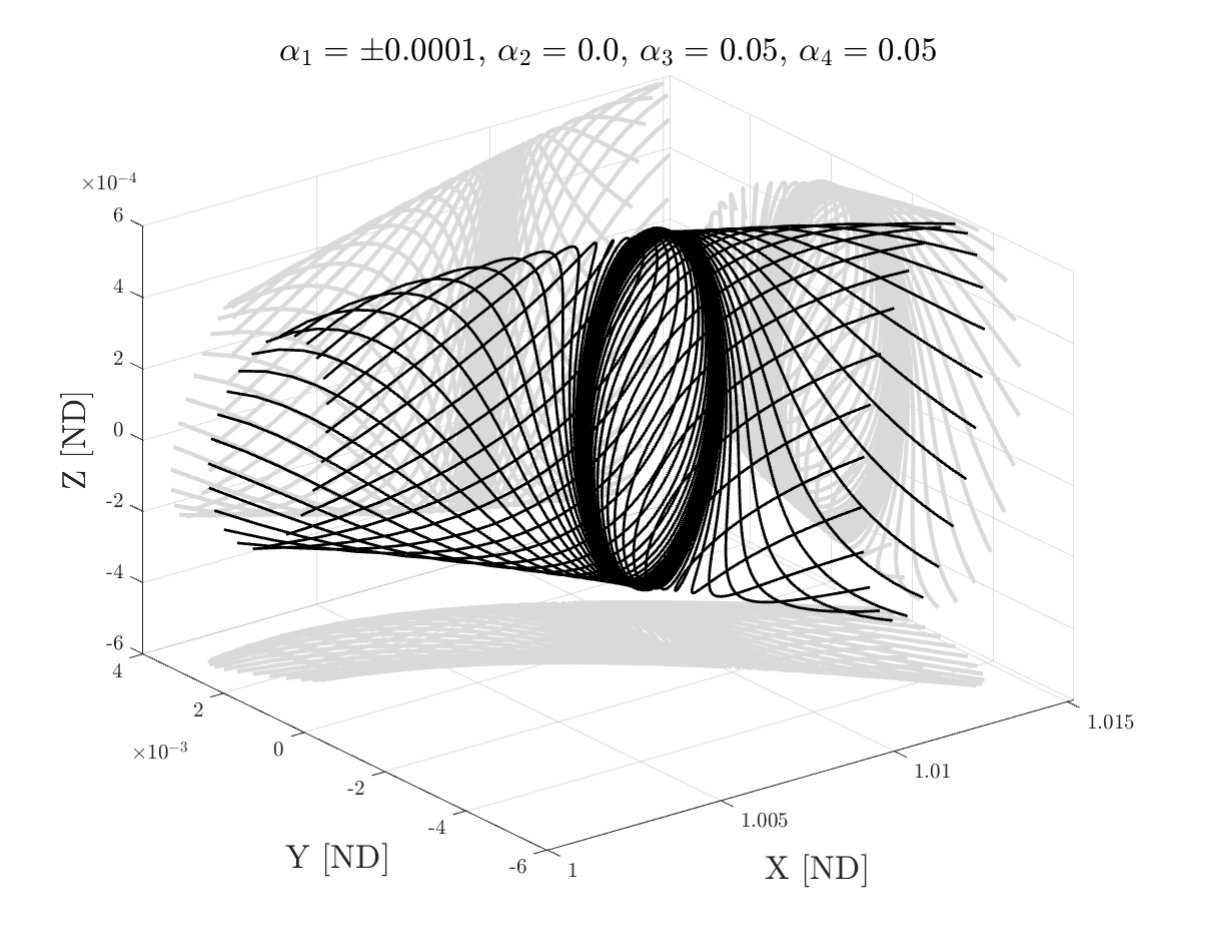}
}
\caption{Manifolds for $\alpha = 80^\circ$, $\gamma = 40^\circ$}
\label{fig:results3}
\end{figure}

In addition to periodic orbits, the approximate solutions may be used to produce stable and unstable manifolds, which lead to and from the AEP, respectively. In Figure \ref{fig:results3}, $\alpha_3 = \alpha_4 = 0.05$ for both plots. In the lefthand plot, $\alpha_1 = 0 $ and $\alpha_2 =  \pm 1 \times 10^{-4}$. In the righthand plot, $\alpha_1 = \pm 1 \times 10^{-4}$ and $\alpha_2 = 0$. Similar to the periodic orbits, these manifolds uniquely exist for \textit{this AEP}, suggesting that a change in one or both of the attitude angles can produce a unique set of manifolds that may better suit a specific mission's purpose.
\begin{figure}[!htp]
\centering
\subfloat[$\alpha = 80^\circ$ and $\gamma = 0^\circ$]{%
  \includegraphics[width=45mm]{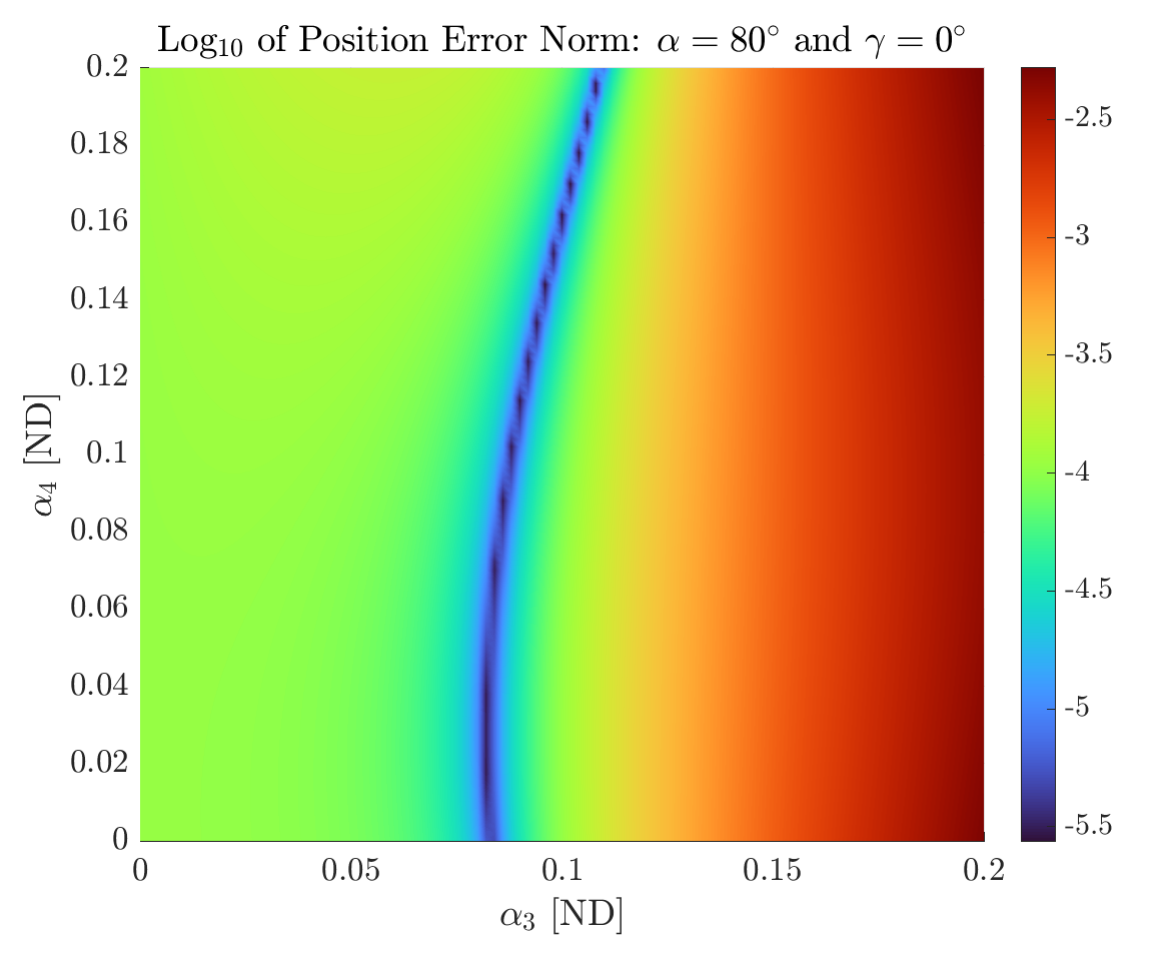}\label{fig:error_alpha80}%
}\quad
\subfloat[$\alpha = 0^\circ$ and $\gamma = 40^\circ$]{%
  \includegraphics[width=45mm]{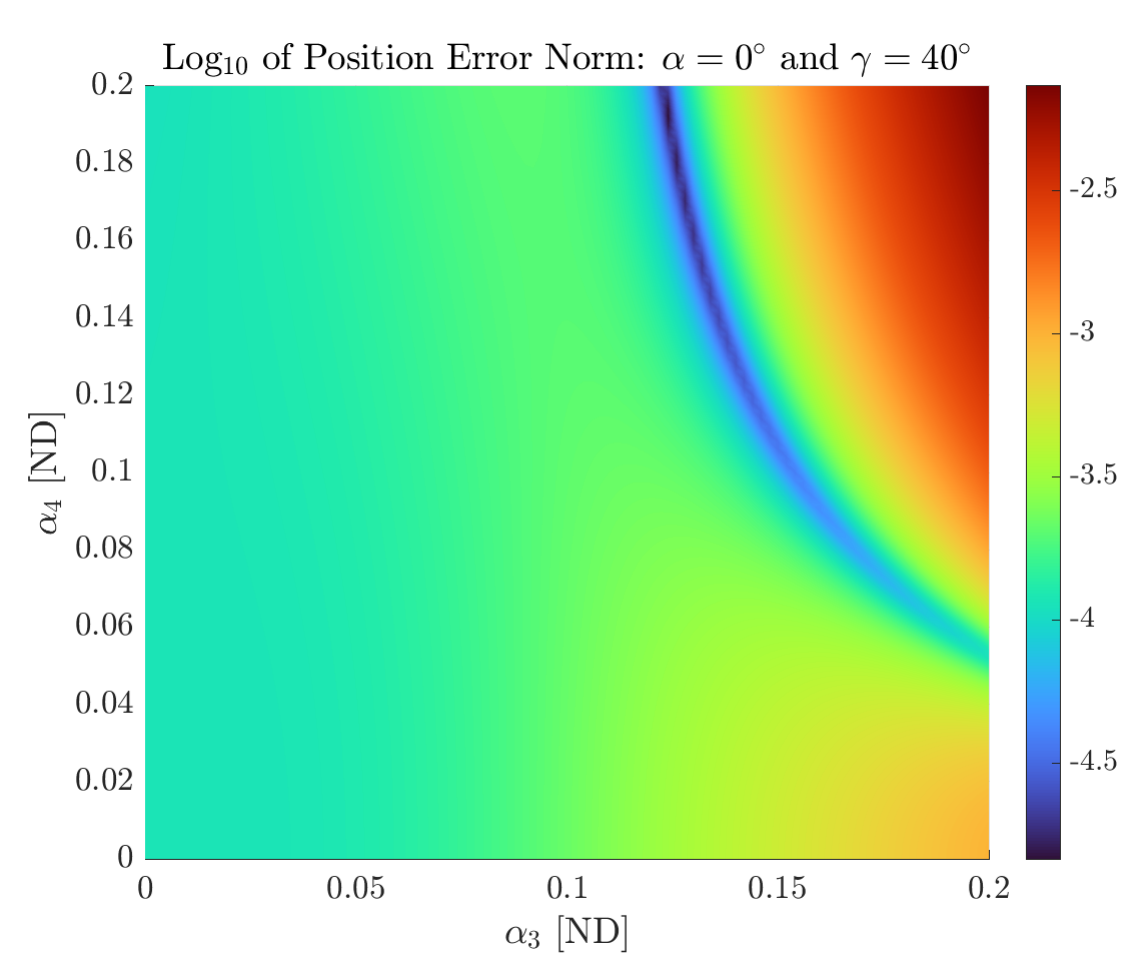}\label{fig:error_gamma40}
}\quad
\subfloat[$\alpha = 80^\circ$ and $\gamma = 40^\circ$]{%
  \includegraphics[width=45mm]{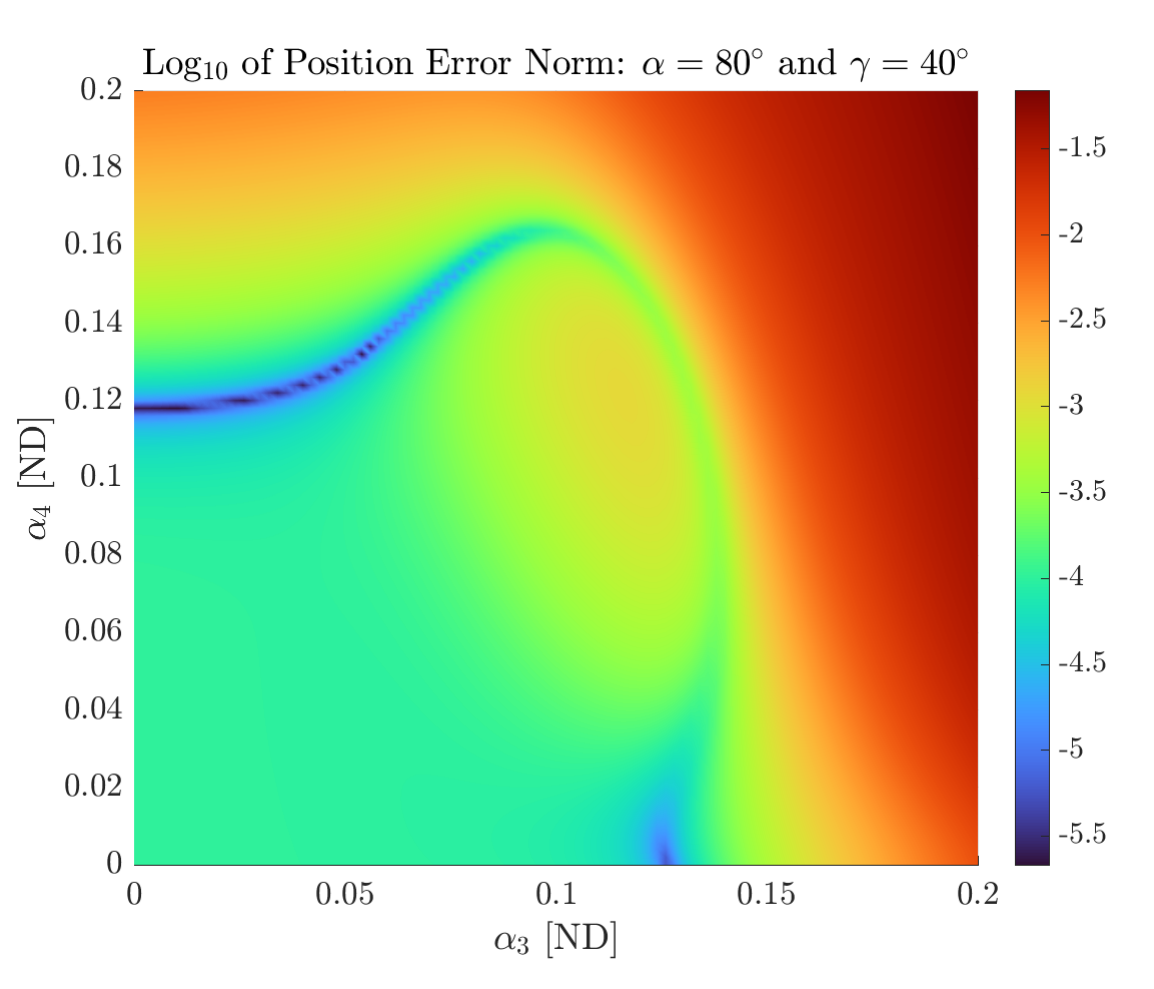}\label{fig:error_alpha80gamma40}
}
\caption{Position Error Norms for Seventh-Order Approximations}
\label{fig:pne}
\end{figure}

\subsection{Position Error Norm}
The accuracy of the seventh-order approximations is determined by using their initial conditions to initialize a numerical integration of the equations of motion (Eq.~\eqref{eqn:eoms}) and comparing the position error norm at $\pi/2$ nondimensional time units. The integration method used was Tsitouras 5/4 Runge-Kutta method in Julia v1.8.5, with absolute and relative error tolerances of $1\times10^{-16}$. The trajectories use $\phi_1 = \phi_2 = 0$ and $\alpha_1 = \alpha_2 = 0$. Additionally, a $100 \times 100$ mesh of $\alpha_3$ and $\alpha_4$ values ranging from $[0,0.2]$ were used to produce trajectories. The Euclidean norm of the errors in position between the numerical solution and the approximate analytical solution is illustrated using logarithmic base 10 for easier visualization in Figure \ref{fig:pne}. It is clear there exist some general patterns between the different sets of attitude angles: 
\begin{enumerate}
	\item The errors generally increase as the amplitudes increase.
	\item There is a curve in which the errors decrease significantly. The relationship between $\alpha_3$ and $\alpha_4$ in this regime is defined and unique for each set of attitude angles. 
	\item Shortly after the regime of low error, the errors begin increasing again.
\end{enumerate}

\section{Conclusion}
This paper detailed the formulation of the CR3BP with SRP equations of motion about an AEP of interest for the non-Hamiltonian SRP case. Furthermore, this paper used the linearized equations of motion to find a first-order approximation of periodic orbits and invariant manifolds. This first-order approximation was then used by the Lindstedt-Poincar\'e method to iteratively determine a user-defined order of approximate analytical, closed-form solutions for both periodic orbits and invariant manifolds. The results shown correspond to seventh-order solutions for $(\alpha,\gamma) = (80^\circ, 0^\circ)$, $(\alpha,\gamma) = (0^\circ, 40^\circ)$, and $(\alpha,\gamma) = (80^\circ, 40^\circ)$ when $\beta = 0.002$ and the AEP of interest is the one analogous to Sun-Earth $L_2$. Both periodic orbits and manifolds were found and illustrated for these sets of attitude angles. Additionally, a numerical study was presented which produced position error norms comparing the approximate analytical solution's position at time $\pi/2$ to that of the numerical, true solution when seeded with the initial condition from the approximation. This study indicates that the error generally increases as the amplitudes increase, but there exist unique regimes of low errors for the sets of attitude angles. Future work will investigate higher-order solutions and a variety of attitude angles and sail lightness numbers so that any restrictions that exist on the utility of the approximate solutions may be well understood and documented. Additionally, future work will incorporate solving for the frequency series for $\omega_r$ and $\nu_r$ rather than assuming they are zero, such that more accurate approximations may be rendered.

\newpage
\bibliographystyle{AAS_publication}
\bibliography{biblio.bib}

\end{document}